\documentclass[10pt]{article}

\usepackage{amsmath,amsxtra,amssymb,latexsym, amscd,amsthm}
\usepackage{here}
\usepackage{graphpap}

\usepackage{graphicx}
\usepackage{color}
\usepackage[]{algorithm2e}
\usepackage{hyperref}

\newtheorem{theorem}{Theorem}[section]
\newtheorem{corollary}[theorem]{Corollary}
\newtheorem{lemma}[theorem]{Lemma}
\newtheorem{proposition}[theorem]{Proposition}
\theoremstyle{definition}

\newtheorem{remark}[theorem]{Remark}

\newtheorem{example}[theorem]{Example}

\numberwithin{equation}{section}

\def\ran{\text{range}}

\def\<{\langle}
\def\>{\rangle}

\def\limn{\lim_{n\to\infty}}
\def\limsupn{\limsup_{n\to\infty}}
\def\liminfn{\liminf_{n\to\infty}}
\definecolor{purple}{rgb}{0.4, 0.0, 0.4}

\begin{document}

\begin{center}
\large \bf A regularization approach for an inverse source problem in elliptic systems
   from single Cauchy data
\end{center}

\centerline { Michael Hinze$^a$\let\thefootnote\relax\footnote{Email:~ michael.hinze@uni-hamburg.de,~ bernd.hofmann@mathematik.tu-chemnitz.de,~ quyen.tran@uni-goettingen.de}, Bernd Hofmann$^b$ and Tran Nhan Tam Quyen$^c$}

{\small
$^a$Department of Mathematics, University of Hamburg, 20146 Hamburg, Germany\\
$^b$Faculty of Mathematics, Chemnitz University of Technology, 09107 Chemnitz, Germany\\
$^c$Institute for Numerical and Applied Mathematics, University of Goettingen, 37083 Goettingen, Germany
}

\centerline{\rule{8cm}{.1pt}}

\vspace{0.3cm}

{\small {\bf Abstract:} In this paper we investigate the problem of identifying the source term $f$ in the elliptic system
$$-\nabla \cdot \big(Q \nabla \Phi \big) = f \mbox{~in~} \Omega \subset R^d, d\in\{2,3\},\, Q \nabla \Phi \cdot \vec{n} = j \mbox{~on~} \partial\Omega\, \mbox{~and~} \Phi = g \mbox{~on~} \partial\Omega$$
from a single noisy measurement couple $\left( j_\delta, g_\delta\right) $ of the Neumann and Dirichlet data $(j,g)$ with noise level $\delta>0$. In this context, the diffusion matrix $Q$ is given. A variational method of Tikhonov-type regularization with specific misfit term of Kohn-Vogelius-type and quadratic stabilizing penalty term is suggested to tackle this linear inverse problem. The method also appears as a variant of the Lavrentiev regularization. For the  occurring linear inverse problem in infinite dimensional Hilbert spaces, convergence and rate results can be found from the general theory of classical Tikhonov
and Lavrentiev regularization. Using the variational discretization concept, where the PDE is discretized with piecewise linear and continuous finite elements, we show the convergence of finite element approximations to solutions
of the regularized problem. Moreover, we derive an error bound and corresponding convergence rates provided a suitable range-type source condition is satisfied. For the numerical solution we propose a conjugate gradient method. To illustrate the theoretical results, a numerical case study is presented which supports our analytical findings.}



\section{Introduction}

Let $\Omega$ be an open, bounded and connected domain of $R^d,
~d\in\{2,3\},$ with Lipschitz boundary $\partial \Omega$. We consider the elliptic system
\begin{align}
-\nabla \cdot \big(Q \nabla \Phi \big) &= f \mbox{~in~} \Omega,  \label{17-5-16ct1}\\
Q \nabla \Phi \cdot \vec{n} &= j^\dag \mbox{~on~} \partial\Omega \mbox{~and~} \label{17-5-16ct2}\\
\Phi &= g^\dag \mbox{~on~} \partial\Omega, \label{17-5-16ct3}
\end{align}
where $\vec{n}$ is the unit outward normal on $\partial\Omega$ and the diffusion matrix $Q$ is given. Furthermore, we assume that $Q :=  \left(q_{rs}\right)_{1\le r, s\le d} \in {L^{\infty}(\Omega)}^{d \times d}$ is symmetric and satisfies the uniformly ellipticity condition
\begin{align}\label{5-5-16ct1}
Q(x)\xi \cdot\xi = \sum_{1\le r,s\le d} q_{rs}(x)\xi_r\xi_s \ge \underline{q} |\xi|^2 \mbox{~a.e. in~} \Omega
\end{align}
for all $\xi = \left(\xi_r\right)_{1\le r\le d} \in R^d$ with some constant $\underline{q} >0$.

The system \eqref{17-5-16ct1}--\eqref{17-5-16ct3} is overdetermined, i.e. if the Neumann and Dirichlet boundary conditions $j^\dag \in H^{-1/2}(\partial\Omega) := {H^{1/2}(\partial\Omega)}^*, ~g^\dag \in H^{1/2}(\partial\Omega)$, and the source term $f \in L^2(\Omega)$ are given, then there may be no $\Phi$ satisfying this system. In this paper we assume that the system is consistent and our aim is to reconstruct a function $f\in L^2(\Omega)$
%and a function $\Phi \in H^1(\Omega)$
in the system \eqref{17-5-16ct1}--\eqref{17-5-16ct3} from a noisy measurement couple $\left(j_\delta,g_\delta \right)\in H^{-1/2}(\partial\Omega) \times H^{1/2}(\partial\Omega) $ of the exact Neumann and Dirichlet data $\big(j^\dag,g^\dag\big)$, where  $\delta >0$ stands for the measurement error, i.e. we assume the noise model
\begin{align}\label{26-3-16ct1}
\big\|j_\delta-j^\dag\big\|_{H^{-1/2}(\partial\Omega)} + \big\|g_\delta-g^\dag\big\|_{H^{1/2}(\partial\Omega)} \le \delta.
\end{align}
The source identification problem in PDEs arises in many branches of applied science such as electroencephalography, geophysical prospecting and pollutant detection, and attracted great attention
from many scientists in the last 30 years or so. For surveys on this subject we may consult in \cite{banks-kunisch1989,chavent2009,EHN96,hama,isakov89,Tarantola} and the references therein. Up to now, only a limited number of works was investigated the general source identification problem and obtained results concentrated on numerical analysis for the identification problem. In \cite{farcas,matsumoto1,matsumoto2} authors have used the dual reciprocity boundary element methods to simulate numerically for the above mentioned identification problem. In case some priori knowledge of the identified source is available, such as a point source, a characteristic function or a harmonic function, numerical methods treating the problem have been obtained in \cite{badia1,badia2,kunisch,ring}. A survey of the problem of simultaneously identifying the source term and coefficients in elliptic systems from {\it distributed} observations can be found in \cite{quyen}, where further references can be found.

In the present paper, the general source identification problem in elliptic partial differential equations from a single noisy measurement couple of Neumann and Dirichlet data is studied. So far, we have not yet found investigations on the discretization analysis for this source recovery problem, a fact which also motivated the research presented in the paper. By using a suitable version of the Tikhonov-type regularization with some non-standard misfit term we could outline that the source distribution inside the physical domain $\Omega$ can be reconstructed from a {\it finite} number of observations on the boundary $\partial\Omega$, at least by numerical approximations. The specific regularization approach proves to be a version of Lavrentiev regularization with implicit forward operator. One of the main results of the paper is to show convergence of the finite element discretized Tikhonov-regularized solutions to a sought source function. Another main result is the interpretation of an occurring condition of solution smoothness as a range-type source condition of Lavrentiev's regularization method. This allows us to establish error bounds and corresponding convergence rates for the regularized solutions.

To formulate precisely the problem, we first give some notations. Let us denote by
$\gamma : H^1(\Omega) \to H^{1/2}(\partial\Omega)$
the continuous Dirichlet trace operator with
$\gamma^{-1} : H^{1/2}(\partial\Omega) \to H^1(\Omega)$
its continuous right inverse operator, i.e. $ (\gamma\circ \gamma^{-1}) g =g$ for all $g\in H^{1/2}(\partial\Omega)$. We set
$$H^1_\diamond(\Omega) := \left\{ u \in H^1(\Omega) ~\Big|~ \int_{\partial\Omega} \gamma udx =0\right\} \mbox{~and~} H^{1/2}_\diamond(\partial\Omega) := \left\{ g \in H^{1/2}(\partial\Omega) ~\Big|~ \int_{\partial\Omega} g(x)dx =0\right\}$$
and denote by $C_\Omega$ the positive constant appearing in the Poincar\'e-Friedrichs inequality (cf.\ \cite{Pechstein})
\begin{align}\label{20-10-15ct1}
C_\Omega \int_\Omega \varphi^2 \le \int_\Omega |\nabla \varphi|^2 \mbox{~for all~} \varphi\in H^1_\diamond(\Omega).
\end{align}
Since $H^1_0(\Omega) := \left\{ u \in H^1(\Omega) ~|~ \gamma u =0\right\} \subset H^1_\diamond(\Omega)$, the inequality \eqref{20-10-15ct1} is in particular valid for all $\varphi\in H^1_0(\Omega)$. Furthermore,
by \eqref{5-5-16ct1}, the coercivity condition
\begin{align}\label{coercivity}
\| \varphi\|^2_{H^1(\Omega)} \le \frac{1+ C_\Omega}{C_\Omega} \int_\Omega | \nabla \varphi|^2 \le \frac{1+ C_\Omega}{C_\Omega \underline{q}} \int_\Omega Q\nabla\varphi \cdot \nabla\varphi
\end{align}
holds for all $\varphi\in H^1_\diamond(\Omega)$.

Now, for any fixed $\left(j,g\right)\in H^{-1/2}(\partial\Omega) \times H^{1/2}_\diamond(\partial\Omega)$ we can simultaneously consider the Neumann problem
\begin{align}
-\nabla \cdot (Q\nabla u) = f \mbox{~in~} \Omega \mbox{~and~}
Q \nabla u \cdot \vec{n} = j \mbox{~on~} \partial\Omega \label{9-6-16ct2}
\end{align}
as well as the Dirichlet problem
\begin{align}
-\nabla \cdot (Q\nabla v) = f \mbox{~in~} \Omega \mbox{~and~} v = g \mbox{~on~} \partial\Omega.  \label{9-6-16ct3}
\end{align}
By the aid of \eqref{coercivity} and the Riesz representation theorem, we conclude that for each
$f\in L^2(\Omega)$
there exists a unique weak solution $u$ of the problem (\ref{9-6-16ct2}) in the sense that $u\in H^1_\diamond(\Omega)$ and satisfies the identity
\begin{align}
\int_\Omega  Q \nabla u \cdot \nabla \varphi =
\left\langle j,\gamma \varphi\right\rangle + \left( f,\varphi\right) \label{ct9}
\end{align}
for all $\varphi\in H^1_\diamond(\Omega)$, where notation $\left\langle j,g\right\rangle$ stands for the value of the function $j\in H^{-1/2}(\partial\Omega)$ at $g\in H^{1/2}(\partial\Omega)$ and the notation $\left( f,\varphi\right)$ is the inner product of $f$ and $\varphi$ in the space $L^2(\Omega)$.
Then we can define the {\it Neumann operator}
\begin{align*}
\mathcal{N} : L^2(\Omega)
\rightarrow H^1_\diamond(\Omega) \mbox{~with~} f\mapsto \mathcal{N}_fj,
\end{align*}
which maps each $f \in L^2(\Omega) $ to the unique weak solution $\mathcal{N}_fj := u$ of the problem (\ref{9-6-16ct2}). Similarly, the problem (\ref{9-6-16ct3}) also attains a unique weak solution $v$ in the sense that $v\in H^1(\Omega)$, $\gamma v = g$ and the identity
\begin{align}
\int_\Omega  Q \nabla v \cdot \nabla \psi  =
\left( f,\psi\right)  \label{ct9*}
\end{align}
holds for all $\psi\in H^1_0(\Omega)$. The {\it Dirichlet operator} is defined as
\begin{align*}
\mathcal{D} : L^2(\Omega)
\rightarrow H^1_\diamond(\Omega) \mbox{~with~} f\mapsto\mathcal{D}_fg,
\end{align*}
which maps each $f \in L^2(\Omega)$ to the unique weak solution $\mathcal{D}_fg := v$ of the problem (\ref{9-6-16ct3}). Therefore, for any fixed $f \in L^2(\Omega)$ we can define the so-called {\it Neumann-to-Dirichlet map}
\begin{align*}
\Lambda_f: H^{-1/2}(\partial\Omega)
\rightarrow H^{1/2}_\diamond(\partial\Omega), \quad
j\mapsto \Lambda_fj := \gamma \mathcal{N}_fj.
\end{align*}

We mention that since $H^1_0(\Omega) \subset H^1_\diamond(\Omega)$, we from \eqref{ct9} have that
$
\int_\Omega  Q \nabla \mathcal{N}_f j \cdot \nabla \psi = \left( f, \psi\right)
$
for all $\psi \in H^1_0(\Omega)$. In view of \eqref{ct9*} we therefore conclude
$\Lambda_fj = g \mbox{~if and only if~} \mathcal{N}_fj = \mathcal{D}_fg$,
where the identities
\begin{align}\label{Plon1}
\mathcal{N}_f j =\mathcal{N}_f0 + \mathcal{N}_0 j \quad \mbox{and} \quad \mathcal{D}_f g =\mathcal{D}_f0 + \mathcal{D}_0 g
\end{align}
are satisfied, and the operators $f\mapsto \mathcal{N}_f0$ and  $f\mapsto \mathcal{D}_f0$ are linear and bounded from $L^2(\Omega)$ into
itself. Furthermore,
$\Lambda_f j = \gamma \mathcal{N}_fj = \gamma\mathcal{N}_0j + \gamma\mathcal{N}_f0 = \Lambda_0j+\Lambda_f0$,
where $\Lambda_0j$ is linear, self-adjoint, bounded and invertible, as the diffusion $Q$ is smooth enough (cf.~\cite{nachman}).

As in electrical impedance tomography (EIT) or for the Calder\'on's problem \cite{astala,calderon,nachman} one can pose the question whether the source distribution $f$ inside a physical domain $\Omega$ can be determined from an {\it infinite} number of observations on the boundary $\partial
\Omega$, i.e. from the Neumann-to-Dirichlet map $\Lambda_f$:
\begin{align*}
f_1, f_2 \in L^2(\Omega) \quad \mbox{with} \quad \Lambda_{f_1} = \Lambda_{f_2}
\quad \Rightarrow \quad f_1=f_2 \ ?
\end{align*}
To the best of our knowledge, the above question is still open so far. In case an observation $\Lambda_\delta$ of $\Lambda_f$ being available one can use a certain regularization method to approximate the sought source. For example, one can consider for operator norms $\|\cdot\|_*$ a minimizer of the problem
$$\min_{f\in L^2(\Omega)} \|\Lambda_f - \Lambda_\delta\|_*^2 + \rho\|f-f^*\|_{L^2(\Omega)}^2$$
as a reconstruction along the lines of Tikhonov's regularization method, where $\rho >0$ is the regularization parameter and $f^*$ is an a-priori estimate of the sought source.

However, in practice we have only a {\it finite} number of observations and the task is to reconstruct the identified source, at least by numerical approximations. Furthermore, for simplicity of exposition we below restrict ourselves to the case of just one observation pair $\left(j_\delta,g_\delta\right)$ being available, while the approach described here can be easily extended to multiple measurements $\left(j_\delta^i,g_\delta^i \right)_{i=1,\ldots,I}$, see Section \ref{Numer_text}, Ex. \ref{ex_mult}. The inverse problem is thus stated as follows.
$$
 Given~ \left( j^\dag, g^\dag \right) \in H^{-1/2}(\partial\Omega) \times H^{1/2}_\diamond(\partial\Omega) ~with~ \Lambda_{f} j^\dag = g^\dag, ~find~ f\in L^2(\Omega). \eqno\left(\mathcal{IP}\right)
$$
In other words, the interested problem is, for given $\left( j^\dag, g^\dag \right) \in H^{-1/2}(\partial\Omega) \times H^{1/2}_\diamond(\partial\Omega)$, to find some $f \in L^2(\Omega) $ and consequently $\Phi \in H^1_\diamond(\Omega)$ such that the system \eqref{17-5-16ct1}--\eqref{17-5-16ct3} is satisfied in the weak sense.
Precisely, we define the general solution set
\begin{equation}\label{inv}
\mathcal{I} \left(j^\dag,g^\dag\right) :=\left\{ f \in L^2(\Omega) ~\big|~ \Lambda_f j^\dag = g^\dag \right\} =\left\{ f \in L^2(\Omega) ~\big|~ \mathcal{N}_fj^\dag = \mathcal{D}_fg^\dag \right\}
\end{equation}
of the inverse problem $\left(\mathcal{IP}\right)$.
The source identification problem as described here is well known to be not uniquely determined from boundary observations (see a counterexample in \cite{alves}), i.e. the set $\mathcal{I} \left(j^\dag,g^\dag\right)$
fails to be a singleton. Since not the Neumann-to-Dirichlet map is given, but only one pair $\left( j^\dag, g^\dag \right)$, the problem is even highly underdetermined. Thus instead we will search for the uniquely determined $f^*$-minimum-norm solution $f^\dagger$, which is the minimizer of the problem
$$
\min_{f \in \mathcal{I} \left(j^\dag,g^\dag\right) }  \|f-f^*\|^2_{L^2(\Omega)}. \eqno\left(\mathcal{IP}-MN\right)
$$
As a consequence of item (iii) of Lemma~\ref{Plon2} below, the set $\mathcal{I} \left(j^\dag,g^\dag\right)$ is non-empty, closed and convex, hence $f^\dag$ is uniquely determined. On the other hand, for all $f\in \mathcal{I} \left(j^\dag,g^\dag\right)$ the equation $\mathcal{N}_f j^\dag = \mathcal{D}_f g^\dag$ is fulfilled. However, we have to solve this equation with noise data $\left(j_\delta,g_\delta\right)\in H^{-1/2}(\partial\Omega) \times H^{1/2}_\diamond(\partial\Omega)$ of $\left( j^\dag, g^\dag \right)$ satisfying \eqref{26-3-16ct1}. The simplest variety of regularization may be to consider a minimizer of the Tikhonov functional
\begin{align}\label{Malta-1}
\|\mathcal{N}_f j_\delta - \mathcal{D}_f g_\delta\|^2_{L^2(\Omega)} + \rho\|f-f^*\|^2_{L^2(\Omega)}
\end{align}
over $f\in L^2(\Omega)$ as an approximation solution to $f^\dagger$.

In present work we adopt the variational approach of Kohn and Vogelius \cite{Kohn_Vogelius1,Kohn_Vogelius11,Kohn_Vogelius2} in using cost functional containing the gradient of forward operators to the above mentioned inverse source problem. More precisely, we use the convex functional
\begin{align}\label{26-3-16ct2}
\mathcal{J}_\delta(f) := \int_\Omega Q\nabla\left(\mathcal{N}_f j_\delta - \mathcal{D}_f g_\delta\right) \cdot \nabla\left(\mathcal{N}_f j_\delta - \mathcal{D}_f g_\delta\right)dx,
\end{align}
instead of the mapping $f\mapsto\|\mathcal{N}_f j_\delta - \mathcal{D}_f g_\delta\|^2_{L^2(\Omega)}$, together with Tikhonov regularization and consider the {\it unique} solution $f_{\rho,\delta}$ of the strictly convex minimization problem
$$\min_{f \in L^2(\Omega)} \Upsilon_{\rho,\delta}(f) \quad\mbox{with}\quad \Upsilon_{\rho,\delta} (f) := \mathcal{J}_\delta(f) + \rho \|f-f^*\|^2_{L^2(\Omega)}, \eqno \left(\mathcal{P}_{\rho,\delta}\right)$$
where the gradient of the functional $\Upsilon_{\rho,\delta}$ can be explicitly written as
\begin{align}\label{Malta-8}
\frac{1}{2}\nabla \Upsilon_{\rho,\delta} (f) = \mathcal{N}_fj_\delta - \mathcal{D}_fg_\delta + \rho(f-f^*) \quad\mbox{for all}\quad f\in L^2(\Omega).
\end{align}
The motivation in using this cost functional $\mathcal{J}_\delta$ as misfit functional is that for all $\xi \in L^2(\Omega)$ the inequality
\begin{align*}
\mathcal{J}_0(\xi) &:= \int_\Omega Q\nabla\left(\mathcal{N}_\xi j^\dag - \mathcal{D}_\xi g^\dag\right) \cdot \nabla\left(\mathcal{N}_\xi j^\dag - \mathcal{D}_\xi g^\dag\right)dx \ge \frac{C_\Omega \underline{q}}{1+ C_\Omega} \left\|\mathcal{N}_\xi j^\dag - \mathcal{D}_\xi g^\dag\right\|^2_{H^1(\Omega)} \ge 0
\end{align*}
holds true and $\mathcal{J}_0(f) =0$ at any $f\in \mathcal{I} \left(j^\dag,g^\dag\right)$. The advantage is evident, because the minimizer $f_{\rho,\delta} \in L^2(\Omega)$ satisfies the equation $\nabla\Upsilon_{\rho,\delta}(f_{\rho,\delta})=0$
such that, for $j:=j_\delta$, $g:=g_\delta$ and $f:=f_{\rho,\delta}$, we have
\begin{equation}\label{Bernd-1}
\mathbf{L}_{j,g}(f):=\mathcal{N}_fj - \mathcal{D}_fg + \rho(f-f^*)=0,
\end{equation}
and hence, for $f:=f_{\rho,\delta}$, we have
\begin{align}
f -f^* =-\frac{1}{\rho} \left( \mathcal{N}_fj_\delta - \mathcal{D}_fg_\delta\right).\label{10-5-16ct3**}
\end{align}
Due to formula (\ref{10-5-16ct3**}), the Tikhonov regularization approach under consideration with specific misfit term also appears as a variant of the Lavrentiev regularization (see, e.g., \cite{AlbRya06,BoHo16,HKR16,Taut02}). After some operator-theoretic settings and preliminary results in Section~2, concerning also the ill-posedness of the linear inverse problem under consideration, we apply in Section~3 the general theory of classical Tikhonov and Lavrentiev regularization for such problems yielding propositions on convergence and convergence rates for the regularized solutions in the infinite dimensional Hilbert spaces. However,
for convenience in numerical analysis with the finite element methods introduced in Section~4 our focus is here on the extremal problem for minimizing the Tikhonov functional with Kohn-Vogelius misfit term and quadratic penalty. The use of different convex penalty terms, e.g. total variation, may be a work for us in future.

Let $\mathcal{N}^h_f j_\delta$ and $\mathcal{D}^h_f g_\delta$ be corresponding approximations of the solution maps $\mathcal{N}_f j_\delta$ and $\mathcal{D}_f g_\delta$ in the finite dimensional space $\mathcal{V}^h_1$ of piecewise linear, continuous finite elements. We then consider the discrete regularized problem corresponding to $\left(\mathcal{P}_{\rho,\delta}\right)$, i.e., the following strictly convex minimization problem
$$
\min_{f \in L^2(\Omega)} \int_\Omega Q\nabla\left(\mathcal{N}^h_f j_\delta - \mathcal{D}^h_f g_\delta\right) \cdot \nabla\left(\mathcal{N}^h_f j_\delta - \mathcal{D}^h_f g_\delta\right)dx + \rho \|f-f^*\|^2_{L^2(\Omega)}. \eqno \left(\mathcal{P}^h_{\rho,\delta}\right)
$$
Using the variational discretization concept introduced in \cite{Hinze}, we show in Section \ref{Stability} that the unique solution $f^h_{\rho,\delta}$ of the problem $\left(\mathcal{P}^h_{\rho,\delta}\right)$ automatically belongs to the finite dimensional space $\mathcal{V}^h_1$. Thus, a discretization of the admissible set $L^2(\Omega)$ can be avoided.

As $h, \delta \to 0$ and with an appropriate a-priori regularization parameter choice $\rho = \rho(h,\delta)$, also in Section \ref{Stability}, we prove that the sequence $\big(f^h_{\rho,\delta}\big)$ converges to $f^\dag$ in the $L^2(\Omega)$-norm. Furthermore, the corresponding state sequences $\Big(\mathcal{N}^h_{f^h_{\rho,\delta}}j_\delta\Big)$ and $\Big(\mathcal{D}^h_{f^h_{\rho,\delta}}g_\delta\Big)$ converge in the $H^1(\Omega)$-norm to $\Phi^\dag =\Phi^\dag(f^\dag,j^\dag,g^\dag)$ solving \eqref{17-5-16ct1}--\eqref{17-5-16ct3}.

Section \ref{tdht} is devoted to convergence rates for the discretized problem. In this section we also show that if $f\in \mathcal{I} \left(j^\dag,g^\dag\right)$ and there is a function $w\in L^2(\Omega)$ such that $f -f^* =\mathbf{L}_{j^\dagger,g^\dagger}(w)$, or in other notation
\begin{align}\label{4-1-17ct1}
f -f^* = \mathcal{N}_wj^\dag - \mathcal{D}_wg^\dag,
\end{align}
then $f=f^\dag$, i.e. $f$ is the {\it unique} $f^*$-minimum-norm solution of the identification problem. Condition (\ref{4-1-17ct1}) appears to be a {\it source condition} for both, Tikhonov and Lavrentiev regularization,
and allows for corresponding convergence rates of the continuous setting in infinite dimensional spaces as well as after incorporating the discretization. In the latter case, precisely for the known matrix $Q\in {C^{0,1}(\Omega)}^{d\times d}$ and the exact data $\left( j^\dag,g^\dag\right)  \in H^{1/2}(\partial\Omega) \times H^{3/2}(\partial\Omega)$, we derive the convergence rates
$$\left\| \mathcal{N}^h_{f^h_{\rho,\delta}} j_\delta - \mathcal{D}^h_{f^h_{\rho,\delta}} g_\delta\right\|^2_{H^1(\Omega)}  + \rho\big\|f^h_{\rho,\delta} - f^\dag \big\|^2_{L^2(\Omega)} =\mathcal{O} \left(\delta^2 + h^2 + h\rho + \delta\rho + \rho^2 \right)$$
and
$$\Big\|\mathcal{N}^h_{f^h_{\rho,\delta}}j_\delta - \Phi^\dag \Big\|^2_{H^1(\Omega)} + \left\|\mathcal{D}^h_{f^h_{\rho,\delta}}g_\delta - \Phi^\dag \right\|^2_{H^1(\Omega)} =\mathcal{O} \left(\delta^2\rho^{-1} + h^2\rho^{-1} + h + \delta+\rho \right).$$
Finally, for the numerical solution of the discrete regularized problem $\left(\mathcal{P}^h_{\rho,\delta}\right)$ we employ in Section~6 a conjugate gradient algorithm. Numerical case studies illustrate the analytical results and show the efficiency of our theoretical findings.

We conclude this introduction with a remark that since the main interest is to clearly state our ideas, we only treat the model elliptic problem \eqref{17-5-16ct1} while the approach described here can be easily extended to more general models, e.g., for the source identification problem in diffusion-reaction equations
\begin{align}\label{27-3-17ct1}
-\nabla \cdot \big(Q \nabla \Phi \big) + \kappa^2 \Phi = f \mbox{~in~} \Omega, ~Q\nabla \Phi\cdot \vec{n} + \sigma\Phi= j^\dag \mbox{~on~} \partial\Omega \mbox{~and~} \Phi = g^\dag \mbox{~on~} \partial\Omega
\end{align}
from a measurement $\left( j_\delta, g_\delta\right) $ of $\big(j^\dag, g^\dag\big)$, where $Q$ satisfying the condition \eqref{5-5-16ct1}, $0 \neq \kappa = \kappa(x) \in L^\infty(\Omega)$, i.e the set $\{x \in \Omega | \kappa(x) \neq 0\}$ has positive Lebesgue measure, and $\sigma =\sigma(x) \in L^\infty(\partial\Omega)$ with $\sigma\ge 0$ are given. The variational approach is now formulated as the minimizing problem with the misfit
$$\int_\Omega Q\nabla\left(R_f j_\delta - D_f g_\delta\right) \cdot \nabla\left(R_f j_\delta - D_f g_\delta\right)dx + \int_\Omega \kappa^2 \left(R_f j_\delta - D_f g_\delta\right)^2dx +\int_{\partial\Omega}\sigma \left(R_f j_\delta - D_f g_\delta\right)^2dx$$
over $f \in L^2(\Omega)$, where $R$ and $D$ are the Robin operator and the Dirichlet operator relating with the equation \eqref{27-3-17ct1}, respectively.

We here would like to mention an inverse problem related closely to the identification in this paper, the problem of identifying the source term $f$ in the Helmholtz-type equation
\begin{align*}
\nabla \cdot \big(Q \nabla \Phi \big) + \kappa^2 \Phi = f \mbox{~in~} \Omega
\end{align*}
from measured Cauchy data $\left( j_\delta(\kappa), g_\delta(\kappa)\right) $ which is available for all frequency $\kappa>0$. The uniqueness results for this identification problem can be found in \cite{alves,bao}, while several effective recovered algorithms have been presented in \cite{acosta,batoul}.

Throughout the paper we use the
standard notion of Sobolev spaces $H^1(\Omega)$, $H^1_0(\Omega)$, $W^{k,p}(\Omega)$, etc from, for example, \cite{Troianiello}. If not stated otherwise we write
$\int_\Omega \cdots$ instead of $\int_\Omega \cdots dx$.

\section{Preliminaries and operator-theoretic settings}

In order to define appropriate operators, we recall the decompositions \eqref{Plon1} of the corresponding Neumann and Dirichlet problems, where $\mathcal{N}_f 0$ and $\mathcal{D}_f 0$ characterize linear mappings of $f \in L^2(\Omega)$. On the other hand,
$\mathcal{N}_0 j$ and $\mathcal{D}_0 g$ depend nonlinearly on $j$ and $g$, respectively, but both are independent of $f$. Hence, the difference of $\mathcal{N}_fj$ and $\mathcal{D}_fg$ characterizes, for fixed elements $j$ and $g$, the affine mapping $\mathbf{L}_{j,g}$  of $f \in L^2(\Omega)$ defined by formula (\ref{Bernd-1}). First we introduce the linear operator $\mathbf{\tilde T} : L^2(\Omega) \to H^1_\diamond(\Omega)$ defined as
$$\mathbf{\tilde T}(f) := \mathcal{N}_f0 - \mathcal{D}_f0 \in  H^1_\diamond(\Omega)\,.$$
Since the image elements $\mathcal{N}_f0 - \mathcal{D}_f0 \in H^1_\diamond(\Omega)$ also belong to $L^2(\Omega)$, one can moreover introduce
 the operator $\mathbf{T} : L^2(\Omega) \to L^2(\Omega)$ defined by
\begin{align}\label{Malta-3}
\mathbf{T}(f) :=\mathcal{N}_f0 - \mathcal{D}_f0 \in L^2(\Omega)\,,
\end{align}
where $\mathbf{T}(f)=\mathbf{\tilde T}(f)$ for all $f \in L^2(\Omega)$.
On the other hand, we remark that
the expression
\begin{align}\label{5-6-16ct1}
[u,v] := \int_\Omega Q\nabla u\cdot \nabla v
\end{align}
generates an inner product on the space $H^1_\diamond(\Omega)$ which is equivalent to the usual one.
Now let $$\mathbf{\tilde T}^* : H^1_\diamond(\Omega) \to L^2(\Omega)$$
be the adjoint operator of $\mathbf{\tilde T} : L^2(\Omega) \to H^1_\diamond(\Omega)$, where $H^1_\diamond(\Omega)$ is equipped with the  inner product \eqref{5-6-16ct1} above. For all $f\in L^2(\Omega)$ and $\phi\in H^1_\diamond(\Omega)$ we thus have
\begin{align}\label{17-8-17ct1}
\left[ \mathbf{\tilde T}f, \phi\right] = \int_\Omega Q\nabla\mathcal{N}_f0\cdot\nabla\phi - \int_\Omega Q\nabla\mathcal{D}_f0\cdot\nabla\phi = (f,\phi) - \int_\Omega Q\nabla\mathcal{D}_f0\cdot\nabla\phi = (f, \mathbf{\tilde T}^*\phi),
\end{align}
by \eqref{ct9}. We now decompose $H^1_\diamond(\Omega)$ into the orthogonal direct sum $H^1_\diamond(\Omega) = H^1_0(\Omega) \oplus {H^1_0(\Omega)}^\perp$ with respect to the inner product \eqref{5-6-16ct1}. We note for all $g\in H^{1/2}_\diamond(\partial\Omega)$ that $\mathcal{D}_0g \in {H^1_0(\Omega)}^\perp$. Furthermore,
$$\forall~ g_1, ~g_2 \in H^{1/2}_\diamond(\partial\Omega),~ g_1\neq g_2 \quad \Rightarrow \quad \mathcal{D}_0g_1 \neq \mathcal{D}_0g_2$$
which implies $\dim {H^1_0(\Omega)}^\perp \ge \dim H^{1/2}_\diamond(\partial\Omega) = \infty$. For all $f\in L^2(\Omega)$ we deduce from \eqref{ct9*} and \eqref{17-8-17ct1} that
$$\phi\in H^1_0(\Omega) \Leftrightarrow  \int_\Omega Q\nabla\mathcal{D}_f0\cdot\nabla\phi = (f,\phi) \Leftrightarrow (f, \mathbf{\tilde T}^*\phi)=0 \Leftrightarrow \mathbf{\tilde T}^*\phi =0 \Leftrightarrow \phi\in \ker \mathbf{\tilde T}^*,$$
or in other words $\ker \mathbf{\tilde T}^* = H^1_0(\Omega)$. Furthermore, for all $\widehat{\phi}\in {H^1_0(\Omega)}^\perp$ we get $\int_\Omega Q\nabla\mathcal{D}_f0\cdot\nabla\widehat{\phi}=0$, since $\mathcal{D}_f0 \in H^1_0(\Omega)$. Again, the equation \eqref{17-8-17ct1} implies that
$$(f,\widehat{\phi}) = (f, \mathbf{\tilde T}^*\widehat{\phi}) \quad \mbox{for all} \quad f\in L^2(\Omega), ~\widehat{\phi}\in {H^1_0(\Omega)}^\perp.$$
Therefore, $\mathbf{\tilde T}^*_{|{H^1_0(\Omega)}^\perp}$ is the compact embedding ${H^1_0(\Omega)}^\perp \hookrightarrow L^2(\Omega)$ and $\mathbf{\tilde T}^*$ is the composition of the projector from  $H^1_\diamond(\Omega)$ onto ${H^1_0(\Omega)}^\perp$ and the embedding operator from $H^1(\Omega)$ to
$L^2(\Omega)$. Furthermore we have, for all $f \in L^2(\Omega)$,
\begin{align}\label{Malta-30}
T(f) =\mathbf{\tilde T}^*[\mathbf{\tilde T}(f)]\,,
\end{align}
because ${\rm range}(\mathbf{\tilde T})$ is orthogonal to $\ker \mathbf{\tilde T}^*$ and hence $\mathbf{\tilde T}^*$ acts only as embedding operator.

\begin{lemma}\label{Plon2}
(i) The operator $\mathbf{T}$ defined by formula (\ref{Malta-3})
is linear, bounded, self-adjoint and non-negative, i.e.~we have
\begin{align}\label{16-5-18f1}
\left(\mathbf{T}(f),w\right) = \left(f,\mathbf{T}(w)\right) \quad \mbox{and} \quad \left(\mathbf{T}(f),f\right)\ge 0 \quad  \mbox{for all} \quad f,w \in L^2(\Omega).
\end{align}
Moreover, $\mathbf{T}$ is compact and has an infinite dimensional range which is non-closed, i.e.~we have
$\overline{\ran(\mathbf{T})} \not= \ran(\mathbf{T})$.

(ii) For any fixed $(j,g)\in H^{-1/2}(\partial\Omega) \times H^{1/2}(\partial\Omega)$ the map $\mathbf{L}_{j,g} : L^2(\Omega) \to L^2(\Omega)$ defined by
$$\mathbf{L}_{j,g}(f) := \mathcal{N}_fj - \mathcal{D}_fg = \mathbf{T}(f) + \mathcal{N}_0j-\mathcal{D}_0g$$
is affine linear, continuous and monotone, i.e.~we have
\begin{align}\label{16-5-18f2}
\left(\mathbf{L}_{j,g}(f)-\mathbf{L}_{j,g}(w),f-w\right) \ge 0 \quad \mbox{for all} \quad f,w \in L^2(\Omega).
\end{align}

(iii) The solution set $\mathcal{I}$ (cf.~(\ref{inv})) is a closed affine subspace of the Hilbert space $L^2(\Omega)$.
\end{lemma}

\begin{proof}
(i) It follows from \eqref{ct9} that
\begin{align}\label{plon3}
\left(\mathbf{T}(f),w\right) = \int_\Omega Q\nabla\mathcal{N}_w0\cdot\nabla (\mathcal{N}_f0 -\mathcal{D}_f0) = \int_\Omega Q\nabla\mathcal{N}_w0\cdot\nabla\mathcal{N}_f0 - \int_\Omega Q\nabla\mathcal{N}_w0\cdot\nabla \mathcal{D}_f0,
\end{align}
and similarly,\begin{align}\label{plon4}
\left(f,\mathbf{T}(w)\right) = \int_\Omega Q\nabla\mathcal{N}_f0\cdot\nabla\mathcal{N}_w0 - \int_\Omega Q\nabla\mathcal{N}_f0\cdot\nabla \mathcal{D}_w0.
\end{align}
Using \eqref{ct9}--\eqref{ct9*} again, we get
\begin{equation}\label{plon5}
\begin{aligned}
& \int_\Omega Q\nabla\mathcal{N}_w0\cdot\nabla \mathcal{D}_f0 = (w,\mathcal{D}_f0)=\int_\Omega Q\nabla\mathcal{D}_w0\cdot\nabla \mathcal{D}_f0\\
& \int_\Omega Q\nabla\mathcal{N}_f0\cdot\nabla \mathcal{D}_w0 = (f,\mathcal{D}_w0)=\int_\Omega Q\nabla\mathcal{D}_f0\cdot\nabla \mathcal{D}_w0
\end{aligned}
\end{equation}
and the self-adjoint property of $\mathbf{T}$ now follows directly from \eqref{plon3}-\eqref{plon5}. We further have from \eqref{ct9}--\eqref{ct9*} for all $f\in L^2(\Omega)$ that
$
\int_\Omega Q\nabla\mathcal{D}_f0\cdot\nabla (\mathcal{N}_f0 -\mathcal{D}_f0) =0.
$
Combining this with the identity $\left(\mathbf{T}(f),f\right) = \int_\Omega Q\nabla\mathcal{N}_f0\cdot\nabla (\mathcal{N}_f0 -\mathcal{D}_f0)$ we arrive at
$\left(\mathbf{T}(f),f\right) = \int_\Omega Q\nabla(\mathcal{N}_f0 -\mathcal{D}_f0)\cdot\nabla (\mathcal{N}_f0 -\mathcal{D}_f0) \ge 0$.

We now show that $\mathbf{T}$ is compact and has an infinite dimensional range. The operator $\mathbf{T} : L^2(\Omega) \to L^2(\Omega)$ as a composition of a bounded linear operator $\mathbf{\tilde T}$ and a compact embedding operator is compact.
Next, we show that $\dim(\mathbf{\tilde T}) = \infty$. For deriving a contradiction we assume that $\dim \ran(\mathbf{\tilde T}) <\infty$. Then we can write $H^1_\diamond(\Omega) = \ran(\mathbf{\tilde T}) \oplus {\ran(\mathbf{\tilde T})}^\perp$ with respect to the inner product \eqref{5-6-16ct1}. By \eqref{17-8-17ct1}, for all $\varphi\in {\ran(\mathbf{\tilde T})}^\perp$ we get
$\int_\Omega Q\nabla\mathcal{D}_f0\cdot\nabla\varphi = (f,\varphi)$ holding for all $f\in L^2(\Omega)$ which implies that $\varphi\in H^1_0(\Omega)$. Therefore, ${H^1_0(\Omega)}^\perp \subset {\big({\ran(\mathbf{\tilde T})}^\perp\big)}^\perp = \overline{\ran(\mathbf{\tilde T})} = \ran(\mathbf{\tilde T})$ and this yields the contradiction $\infty = \dim {H^1_0(\Omega)}^\perp \le \dim \ran(\mathbf{\tilde T}) < \infty$. Consequently, $\mathbf{T} : L^2(\Omega) \to L^2(\Omega)$ is compact operator and possesses an infinite dimensional range.

(ii) The inequality \eqref{16-5-18f2} follows directly from \eqref{16-5-18f1}.

(iii) Since $\mathbf{T}$ is a bounded linear operator, the solution set can be written as
$$\mathcal{I}(j^\dagger,g^\dagger) =\left\{ f= f_0 \oplus f_0^\perp \in L^2(\Omega) ~\big|~ f_0 \in \ker(\mathbf{T}),\; f_0^\perp=\mathbf{T}^\dagger(\mathcal{D}_0 g^\dagger -\mathcal{N}_0 j^\dagger)\right\}$$
with the Moore-Penrose pseudoinverse $\mathbf{T}^\dagger$ of $\mathbf{T}$. Then the nullspace $\ker(\mathbf{T})$ is a closed subspace and $\mathbf{T}^\dagger(\mathcal{D}_0 g^\dagger -\mathcal{N}_0 j^\dagger)$ a well-defined element in
$L^2(\Omega)$. Consequently, $\mathcal{I}(j^\dagger,g^\dagger)$ is a non-empty closed affine subspace and hence also a convex set in the Hilbert space  $L^2(\Omega)$.
\end{proof}

Due to item (ii) of Lemma~\ref{Plon2}, we can reformulate the identification problem $\left(\mathcal{IP}\right)$ as an operator equation with linear bounded self-adjoint non-negative operator $\mathbf{T}$ mapping in $L^2(\Omega)$. Finding an element $f \in \mathcal{I} \left(j^\dag,g^\dag\right)$ is then equivalent to solving the linear operator equation
\begin{equation} \label{opeq}
\mathbf{T}(f)=\kappa(j^\dagger,g^\dagger), \quad \mbox{where} \quad \kappa(j,g):=\mathcal{D}_0 g -\mathcal{N}_0 j.
\end{equation}
This makes the inverse problem explicit, but we have to take into account that instead of the exact right-hand side $\kappa(j^\dagger,g^\dagger)$ only noisy data $\kappa(j_\delta,g_\delta)$ satisfying (\ref{26-3-16ct1}) for
$j_\delta$ and $g_\delta$ are available. As a consequence of item (i) of Lemma~\ref{Plon2} we see that the equation (\ref{opeq}) formulated in the Hilbert space $L^2(\Omega)$ is {\it ill-posed of type~II} in the sense of Nashed (cf.~\cite{Nashed86}). Stable approximate (regularized) solutions $f_{\rho,\delta}$ to equation (\ref{opeq}) satisfy with $f=f_{\rho,\delta}$ the
auxiliary linear operator equation
\begin{equation} \label{Lavrep}
\mathbf{T}(f) + \rho\,(f-f^*)=\kappa(j_\delta,g_\delta)
\end{equation}
in $L^2(\Omega)$ with some regularization parameter $\rho>0$.

At this point we should recall and note that we have, for all $j_\delta$, $g_\delta$ under consideration and $f \in L^2(\Omega)$, and due to
(\ref{ct9}) -- (\ref{ct9*})
$$\mathbf{T}(f) \in {H^1_0(\Omega)}^\perp \subset H^1_\diamond(\Omega) \qquad \mbox{as well as} \qquad  \kappa(j_\delta,g_\delta) \in {H^1_0(\Omega)}^\perp,$$
which means that the elements $\kappa(j_\delta,g_\delta)$ and $\mathbf{\tilde T}^*[\kappa(j_\delta,g_\delta)]$ in $L^2(\Omega)$ coincide.
Nevertheless, we have to distinguish the two cases $\kappa(j_\delta,g_\delta) \in  H^1_\diamond(\Omega)$  and $\kappa(j_\delta,g_\delta) \in L^2(\Omega)$ in the following lemma.

\begin{lemma}\label{New1}
Under the noise model (\ref{26-3-16ct1}) there is a constant $0 \le \tilde K <\infty$ independent of $\delta$ such that
$$\|\kappa(j_\delta,g_\delta)-\kappa(j^\dagger,g^\dagger)\|_{H^1(\Omega)} \le \tilde K\,\delta.$$
Moreover, there is also a  constant $0 \le  K <\infty$ independent of $\delta$ such that
$$\|\kappa(j_\delta,g_\delta)-\kappa(j^\dagger,g^\dagger)\|_{L^2(\Omega)} \le K\,\delta.$$
\end{lemma}

\begin{proof}
By Lemma~\ref{6-5-16ct5} below the existence of such constant $\tilde K$ in the first estimate of this lemma follows with the settings $h:=0$, $f_1 = f_2 := 0$ and $\tilde K:=\max\{\mathcal{C_N},\mathcal{C_D}\}$. Then we have under the noise model (\ref{26-3-16ct1})
\begin{align*}
\|\kappa(j_\delta,g_\delta)-\kappa(j^\dagger,g^\dagger)\|_{H^1(\Omega)}
&\le \| \mathcal{N}_0{j_\delta} - \mathcal{N}_0{j^\dag}\|_{H^1(\Omega)} + \| \mathcal{D}_0{j_\delta} - \mathcal{D}_0{j^\dag}\|_{H^1(\Omega)}\\
&\le \tilde K \left(\|j_\delta-j^\dag \|_{H^{-1/2}(\partial\Omega)} +  \|g_\delta-g^\dag \|_{H^{1/2}(\partial\Omega)} \right) \le \tilde K \delta.
\end{align*}
By setting $K:=\tilde K \|\mathbf{\tilde T}\|$ the second estimate of the lemma gets established.
This completes the proof.
\end{proof}

We conclude this section by mentioning that the functional $\mathcal{J}_\delta$ defined by \eqref{26-3-16ct2} is convex and weakly sequentially lower semi-continuous. In fact,
the above defined operator $\mathbf{\tilde T} : L^2(\Omega) \to H^1_\diamond(\Omega) $ is compact (see the proof of Lemma~\ref{Plon2}). Using the equivalent inner product \eqref{5-6-16ct1},  we therefore conclude that
\begin{align} \label{Bernd-2}
\mathcal{J}_\delta\big(f\big) = \left[ \mathcal{N}_f j_\delta - \mathcal{D}_f g_\delta, \mathcal{N}_f j_\delta - \mathcal{D}_f g_\delta\right] = \left[ \mathbf{\tilde T}(f) -\kappa(j_\delta,g_\delta), \mathbf{\tilde T}(f)- \kappa(j_\delta,g_\delta)\right],
\end{align}
which shows the convexity and weak sequential lower semi-continuity of the functional $\mathcal{J}_\delta$ and is the basis for a classical Tikhonov regularization approach in the subsequent section. Moreover, we have
as a basis for Lavrentiev regularization in the subsequent section
\begin{align} \label{Bernd-3}
\mathcal{J}_\delta\big(f\big)+\rho\|f-f^*\|^2_{L^2(\Omega)}=\left(\mathbf{T}(f),f \right)+\rho\,\|f\|_{L^2(\Omega)}^2-2\left(f,\kappa(j_\delta,g_\delta)+\rho\,f^*\right)+ {\rm const.},
\end{align}
where the constant is independent of $f$,
and it is well-known that, due to the properties of $\mathbf{T}$ from Lemma~\ref{Plon2}, the unique minimizer $f_{\rho,\delta}$ of this functional coincides with the unique solution of the operator equation (\ref{Lavrep}).

We close
 this section by the following note. As discussed in the Introduction section, instead of using Kohn and Vogelius' function \eqref{26-3-16ct2} we can use the least squares function (cf.\ \eqref{Malta-1}), and then the jointed minimization problem reads as
\begin{align}\label{Malta-2}
\min_{f\in L^2(\Omega)} \Theta(f) \quad\mbox{with}\quad \Theta(f) :=  \|\mathbf{T}(f) - \kappa(j_\delta,g_\delta)\|^2_{L^2(\Omega)} + \rho\|f-f^*\|^2_{L^2(\Omega)}
\end{align}
where $\mathbf{T}(f)$ and $\kappa(j_\delta,g_\delta)$ were defined by \eqref{Malta-3} and \eqref{opeq}, respectively. One can show easily that the problem attains a unique solution $\bar{f}$. Furthermore, for computation this minimizer $\bar{f}$ in practice one must derive the $L^2$-gradient $\nabla \Theta(\bar{f})$. Let $\xi\in L^2(\Omega)$ be arbitrary. We will compute briefly the differential $\Theta'(\bar{f})\xi$ as follows.

We have that $\frac{1}{2} \Theta'(\bar{f})\xi = \left(\mathbf{T}(\bar{f}) - \kappa(j_\delta,g_\delta), \mathbf{T}(\xi)\right) + \rho(\bar{f} - f^*, \xi)$. Consider the adjoint problem
\begin{align}
-\nabla \cdot \big(Q \nabla \bar{\Phi} \big) &= \mathbf{T}(\bar{f}) - \kappa(j_\delta,g_\delta) \mbox{~in~} \Omega,  \label{Malta-4}\\
Q \nabla \bar{\Phi} \cdot \vec{n} &= 0 \mbox{~on~} \partial\Omega.\label{Malta-5}
\end{align}
(We here do not use the homogeneous Dirichlet boundary condition $\bar{\Phi} =0$ on $\partial\Omega$ instead of \eqref{Malta-5}, because in general $\mathbf{T}(\xi) \notin H^1_0(\Omega)$.) Next, we decompose
\begin{align}\label{Malta-6}
\bar{\Phi} = \bar{\Phi}^1 \oplus \bar{\Phi}^2 \in  H^1_0(\Omega) \oplus {H^1_0(\Omega)}^\perp
\end{align}
with respect to the inner product \eqref{5-6-16ct1}. Then we obtain that
\begin{align*}
\left(\mathbf{T}(\bar{f}) - \kappa(j_\delta,g_\delta), \mathbf{T}(\xi)\right)
&= \int_\Omega Q \nabla \bar{\Phi} \cdot \nabla \mathbf{T}(\xi) =  \int_\Omega Q \nabla \mathcal{N}_{\xi}0 \cdot \nabla \bar{\Phi} - \int_\Omega Q \nabla \mathcal{D}_{\xi}0 \cdot \nabla \bar{\Phi}^1 - \underbrace{\int_\Omega Q \nabla \mathcal{D}_{\xi}0 \cdot \nabla \bar{\Phi}^2 }_{=0}\\
&= (\xi, \bar{\Phi}) - (\xi, \bar{\Phi}^1) = (\bar{\Phi}^2, \xi).
\end{align*}
As a result, we arrive at $\frac{1}{2}\nabla\Theta(\bar{f}) = \bar{\Phi}^2 + \rho(\bar{f} - f^*)$. Compared with \eqref{Malta-8}, by utilizing Kohn and Vogelius' function \eqref{26-3-16ct2}, we here avoid any computations for the adjoint problem (\ref{Malta-4})--(\ref{Malta-5}). Furthermore, we avoid computing numerically for the terms $\bar{\Phi}^1$ and $\bar{\Phi}^2$ in the decomposition (\ref{Malta-6}) which seems to be still very difficult for us.

\section{Lavrentiev regularization versus Tikhonov regularization for the continuous problem in infinite dimensional Hilbert spaces} \label{sec:3}

In this section, we consider the error analysis of stable approximate solutions to the continuous problem (\ref{opeq}) in infinite dimensional Hilbert spaces by distinguishing the situations that the right-hand side $\kappa$ is considered as an element in $L^2(\Omega)$ and alternatively as an element in $H^1(\Omega)$.

\smallskip

First our focus is on the Lavrentiev regularization approach based on formula (\ref{Bernd-3}), in which $\kappa(j_\delta,g_\delta) \in L^2(\Omega)$.
The general theory of linear Lavrentiev regularization (see, e.g., \cite{Taut02} and also \cite{AlbRya06,BoHo16,HKR16,Plato16,PMH18}) yields convergence and convergence rates results for the error of regularized solutions $f_{\rho,\delta}$ with respect to the uniquely determined $f^*$-minimizing solution $f^\dagger$ to problem ($\mathcal{IP}-MN$). Taking into account that Lemma~\ref{New1} holds, we immediately derive (see, e.g., \cite[Rem. 3.3]{BoHo16}) the error estimate
\begin{equation} \label{error1}
\|f_{\rho,\delta}-f^\dagger\|_{L^2(\Omega)} \le \rho \|(\mathbf{T}+\rho \mathbf{I})^{-1}(f^\dagger-f^*)\|_{L^2(\Omega)} + \frac{K\delta}{\rho} \le \|f^\dagger-f^*\|_{L^2(\Omega)}+ \frac{K\delta}{\rho}.
\end{equation}
This immediately yields (see \cite[section~2]{HKR16}) the following convergence assertion.
\begin{proposition} \label{proerror1}
For any a priori parameter choice $\rho(\delta)$ of the regularization parameter satisfying the conditions
\begin{equation} \label{error2}
\rho(\delta) \to 0  \qquad \mbox{and} \qquad \frac{\delta}{\rho(\delta)} \to 0 \qquad \mbox{as} \quad \delta \to 0
\end{equation}
we have for a sequence $\delta_n \to 0$, associated data $j_{\delta_n},g_{\delta_n}$, and associated Lavrentiev-regularized solutions $f_{\rho(\delta_n),\delta_n}$ that
 \begin{equation*}
\lim \limits_{n \to \infty} \|f_{\rho(\delta_n),\delta_n}-f^\dagger\|_{L^2(\Omega)}=0,
\end{equation*}
i.e. the regularized solutions are strongly convergent in $L^2(\Omega)$ to the $f^*$-minimum norm solution $f^\dagger$.
\end{proposition}

We can also apply the following well-known result on convergence rates from \cite[Theorem 2.2]{Taut02}:

\begin{proposition} \label{proerror2}
If there is a source element $v \in L^2(\Omega)$ such that the range-type source condition
\begin{equation} \label{error3+}
f^\dagger-f^*=\mathbf{T}(v)
\end{equation}
is satisfied, we have for an a priori choice $\rho(\delta) \sim \sqrt{\delta}$ of the regularization parameter $\rho$ the convergence rate
\begin{equation} \label{error4}
\|f_{\rho(\delta),\delta}-f^\dagger\|_{L^2(\Omega)}= \mathcal{O}(\sqrt{\delta}) \qquad \mbox{as} \qquad \delta \to 0
\end{equation}
of Lavrentiev-regularized solutions.
\end{proposition}

\begin{corollary} \label{proerror3}
The rate result (\ref{error4}) remains true if we have some element $w \in L^2(\Omega)$ such that
\begin{equation} \label{error5}
f^\dag -f^* = \mathcal{N}_wj^\dag - \mathcal{D}_wg^\dag.
\end{equation}
\end{corollary}
\begin{proof}
Evidently we have $\mathcal{N}_wj^\dag - \mathcal{D}_wg^\dag= \mathbf{T}(w)-\kappa(j^\dagger,g^\dagger)$ and $f^\dagger$ is a solution to the operator equation (\ref{opeq}) with exact right-hand side $\kappa(j^\dagger,g^\dagger)$.
Hence (\ref{error5}) can be rewritten as
$f^\dag -f^* = \mathbf{T}(w-f^\dagger)$.
This yields (\ref{error3+}) with the new source element $v:=w-f^\dagger$.
\end{proof}

\begin{remark} \label{rem:Plato}
It was shown by a saturation theorem in \cite{Plato16} that (\ref{error4}) is the best possible rate for linear Lavrentiev regularization, with the exception of singular situations with respect to the forward operator, here $\mathbf{T}$,
and with respect to the solution $f^\dagger$.
\end{remark}

Revisiting the source condition (\ref{error5}), we add the following remarkable result.

\begin{proposition}\label{plon8}
Assume that $f \in \mathcal{I}\left(j^\dag,g^\dag\right)$ and that there is a function $w \in L^2(\Omega)$  such that \linebreak $f -f^* = \mathcal{N}_wj^\dag - \mathcal{D}_wg^\dag$. Then $f$ is the uniquely determined $f^*$-minimum-norm solution of problem $\left(\mathcal{IP}-MN\right)$, i.e.~we have $f =f^\dagger$.
\end{proposition}
\begin{proof}
We have with $\xi \in \left\{ \xi \in L^2(\Omega) ~\big|~ \mathcal{N}_{\xi} j^\dag = \mathcal{D}_{\xi}g^\dag \right\}$ that
\begin{align*}
&\frac{1}{2}\|\xi-f^*\|^2_{L^2(\Omega)} - \frac{1}{2}\|f-f^*\|^2_{L^2(\Omega)} \\
&= \frac{1}{2}\|\xi-f\|^2_{L^2(\Omega)} + (f-f^*,\xi-f) \ge (f-f^*,\xi-f)= \left( \mathcal{N}_wj^\dag - \mathcal{D}_wg^\dag, \xi-f\right) \\
&= \int_\Omega Q \nabla\mathcal{N}_\xi j^\dag\cdot \nabla \left( \mathcal{N}_wj^\dag - \mathcal{D}_wg^\dag\right) - \left\langle j^\dag, \gamma\left( \mathcal{N}_wj^\dag - \mathcal{D}_wg^\dag\right)  \right\rangle \\
&~\quad - \int_\Omega Q \nabla\mathcal{N}_f j^\dag\cdot \nabla \left( \mathcal{N}_wj^\dag - \mathcal{D}_wg^\dag\right) + \left\langle j^\dag, \gamma\left( \mathcal{N}_wj^\dag - \mathcal{D}_wg^\dag\right)  \right\rangle \\
&= \int_\Omega Q \nabla \left( \mathcal{N}_\xi j^\dag - \mathcal{N}_f j^\dag\right) \cdot \nabla \left( \mathcal{N}_wj^\dag - \mathcal{D}_wg^\dag\right).
\end{align*}
Since $\gamma \mathcal{N}_\xi j^\dag = \gamma \mathcal{N}_f j^\dag = g^\dag$, it follows that $\mathcal{N}_\xi j^\dag -  \mathcal{N}_f j^\dag \in H^1_0(\Omega)$. We thus obtain from the last inequality
\begin{align*}
\frac{1}{2}\|\xi-f^*\|^2_{L^2(\Omega)} - \frac{1}{2}\|f-f^*\|^2_{L^2(\Omega)}
&\ge  \int_\Omega Q \nabla \mathcal{N}_wj^\dag \cdot \nabla \left( \mathcal{N}_\xi j^\dag - \mathcal{N}_f j^\dag\right) - \int_\Omega Q \nabla \mathcal{D}_wg^\dag \cdot \nabla \left( \mathcal{N}_\xi j^\dag - \mathcal{N}_f j^\dag\right)\\
&= \left( w, \mathcal{N}_\xi j^\dag - \mathcal{N}_f j^\dag\right) + \left\langle j^\dag, \gamma \left( \mathcal{N}_\xi j^\dag - \mathcal{N}_f j^\dag\right)\right\rangle - \left( w, \mathcal{N}_\xi j^\dag - \mathcal{N}_f j^\dag\right)=0,
\end{align*}
which finishes the proof.
\end{proof}
\begin{remark} \label{rembernd}
Indeed, the statement of Proposition~\ref{plon8} is a special case of the general assertion that source conditions $f-f^*=\mathbf{T}(v),\; v \in X,$ for self-adjoint non-negative bounded linear operators $\mathbf{T}$ in the Hilbert space $X$ can only hold if $f$ is a $f^*$-minimum-norm solution (see \cite[Remark~6]{PMH18}). However, as shown above we have here $\mathcal{N}_wj^\dag - \mathcal{D}_wg^\dag=\mathbf{T}(v)$ with $v=w-f^\dagger \in L^2(\Omega)$.
\end{remark}

Our second alternative focus is on the classical Tikhonov regularization approach based on formula (\ref{Bernd-2}), in which $\kappa(j_\delta,g_\delta) \in H^1(\Omega)$. Then the regularized solution can be established as
\begin{equation} \label{Bernd-4}
f_{\rho,\delta}=(\mathbf{\tilde T}^*\mathbf{\tilde T}+\rho \mathbf{I})^{-1}\left[\mathbf{\tilde T}^*(\kappa(g_\delta,j_\delta))+\rho f^*\right]
\end{equation}
(cf., e.g., \cite[Sec.~5.1]{EHN96}). From (\ref{Bernd-4}) and Lemma~\ref{New1} we easily derive the error estimate
\begin{equation} \label{error10}
\|f_{\rho,\delta}-f^\dagger\|_{L^2(\Omega)} \le \rho \|(\mathbf{\tilde T}^*\mathbf{\tilde T}+\rho \mathbf{I})^{-1}(f^\dagger-f^*)\|_{L^2(\Omega)} + \frac{\tilde K\delta}{2 \sqrt{\rho}} \le \|f^\dagger-f^*\|_{L^2(\Omega)} + \frac{\tilde K\delta}{2 \sqrt{\rho}} ,
\end{equation}
on which the following proposition is based.
\begin{proposition} \label{proerror10}
For any a priori parameter choice $\rho(\delta)$ of the regularization parameter satisfying the conditions
\begin{equation} \label{error2+}
\rho(\delta) \to 0  \qquad \mbox{and} \qquad \frac{\delta}{\sqrt{\rho(\delta)}} \to 0 \qquad \mbox{as} \quad \delta \to 0
\end{equation}
we have for a sequence $\delta_n \to 0$, associated data $j_{\delta_n},g_{\delta_n}$, and associated Tikhonv-regularized solutions $f_{\rho(\delta_n),\delta_n}$ that
 \begin{equation*}
\lim \limits_{n \to \infty} \|f_{\rho(\delta_n),\delta_n}-f^\dagger\|_{L^2(\Omega)}=0,
\end{equation*}
i.e. the regularized solutions are strongly convergent in $L^2(\Omega)$ to the $f^*$-minimum norm solution $f^\dagger$.
\end{proposition}

Furthermore, under the source condition (\ref{error3+}), which is equivalent to
\begin{equation} \label{Bernd-5}
f^\dagger-f^*=[\mathbf{\tilde T}^*\mathbf{\tilde T}](v),
\end{equation}
we find even the bounds
\begin{equation} \label{error11}
\|f_{\rho,\delta}-f^\dagger\|_{L^2(\Omega)} \le \rho \|(\mathbf{\tilde T}^*\mathbf{\tilde T}+\rho \mathbf{I})^{-1}[\mathbf{\tilde T}^*\mathbf{\tilde T}](v)\|_{L^2(\Omega)} + \frac{\tilde K\delta}{2 \sqrt{\rho}} \le \rho
\|v\|_{L^2(\Omega)}+\frac{\tilde K\delta}{2 \sqrt{\rho}},
\end{equation}
which as a consequence of $\mathbf{T}(v)=[\mathbf{\tilde T}^*\mathbf{\tilde T}](v)$ for all $v \in L^2(\Omega)$ immediately yields the rate assertion of the following proposition.

\begin{proposition} \label{proerror20}
If there is a source element $v \in L^2(\Omega)$ such that the range-type source condition (\ref{error3+}) (equivalent to (\ref{error5}) due to Corollary~\ref{proerror3})
is satisfied, then we have for an a priori choice $\rho(\delta) \sim \delta^{2/3}$ of the regularization parameter $\rho$ the convergence rate
\begin{equation} \label{error40}
\|f_{\rho(\delta),\delta}-f^\dagger\|_{L^2(\Omega)}= \mathcal{O}(\delta^{2/3}) \qquad \mbox{as} \qquad \delta \to 0
\end{equation}
of Tikhonov-regularized solutions. If the regularization parameter is chosen as  $\rho(\delta) \sim \delta$, then we obtain also here the convergence rate (\ref{error4}) as in Proposition~\ref{proerror2}.
\end{proposition}
\begin{remark} \label{rem:Groetsch}
It was shown by a saturation theorem of Groetsch 1984 (see, e.g, \cite[Proposition~4.20]{EHN96}) that (\ref{error40}) is the best possible rate for classical linear Tikhonov regularization, with the exception of singular cases. At first glance, it is amazing that the best possible rate (\ref{error40}) in Proposition~\ref{proerror20} is higher than the best possible rate (\ref{error4}) in Proposition~\ref{proerror2}. However, here the classical Tikhonov regularization makes use of the higher smoothness assumption $\kappa(j_\delta,g_\delta) \in H^1(\Omega)$, whereas the version of Lavrentiev regularization
employed here ignores this higher smoothness of the right-hand side of equation (\ref{opeq}) and considers $\kappa(j_\delta,g_\delta)$ only as an element in $L^2(\Omega)$, see the two different noise inequalities in Lemma~\ref{New1}.
\end{remark}

\section{Finite element approximation and convergence for the discretized problem}\label{Stability}

We in Section \ref{sec:3} applied the Lavrentiev regularization and the Tikhonov regularization for the continuous identification problem. In the remain Sections \ref{Stability} -- \ref{Numer_text} we will analyze the problem in finite dimensional spaces. So far we have not yet found investigations on the discretization error in a combination of
both functionals for the fully setting, a fact which motivated the research
presented in this paper.

Let $\left(\mathcal{T}^h\right)_{0<h<1}$ be a family of regular and
quasi-uniform triangulations of the domain $\overline{\Omega}$ with the mesh size $h$.
For the definition of the discretization space of the state
functions let us denote
\begin{equation*}
\mathcal{V}_1^h := \left\{v^h\in C(\overline\Omega)
~|~{v^h}_{|T} \in \mathcal{P}_1(T), ~~\forall
T\in \mathcal{T}^h\right\} \quad \mbox{and} \quad \mathcal{V}_{1,\diamond}^h := \mathcal{V}_1^h \cap H^1_\diamond(\Omega) \mbox{~and~} \mathcal{V}_{1,0}^h := \mathcal{V}_1^h \cap H^1_0(\Omega) \subset \mathcal{V}_{1,\diamond}^h,
\end{equation*}
where $\mathcal{P}_1$ consists all polynomial functions of degree
less than or equal to 1.

\begin{proposition}
(i) Let $f$ be in $L^2(\Omega)$ and $j$ be in $H^{-1/2}(\partial\Omega)$. Then the variational equation
\begin{align}
\int_\Omega Q\nabla u^h \cdot \nabla \varphi^h =  \big( f, \varphi^h \big) + \big\langle j,\gamma\varphi^h\big\rangle \mbox{~for all~} \varphi^h\in
\mathcal{V}_{1,\diamond}^h \label{10/4:ct1}
\end{align}
admits a unique solution $u^h\in \mathcal{V}_{1,\diamond}^h$. Furthermore, the estimate
\begin{align}
\big\|u^h\big\|_{H^1(\Omega)}\le C_\mathcal{N} \left( \left\|f\right\|_{L^2(\Omega)}+
\left\|j\right\|_{H^{-1/2}(\partial\Omega)} \right)  \label{18/5:ct1}
\end{align}
is satisfied. The map $\mathcal{N}^h: L^2(\Omega) \rightarrow
\mathcal{V}_{1,\diamond}^h$ from each $f \in L^2(\Omega)$
to the unique solution $u^h =: \mathcal{N}^h_fj$ of \eqref{10/4:ct1}
is then called the  {\it discrete Neumann operator}.

(ii) Let $f$ be in $L^2(\Omega)$ and $g$ be in $H^{1/2}_\diamond(\partial\Omega)$. The equation
\begin{align}
\int_\Omega Q\nabla v^h \cdot \nabla \psi^h = \big( f, \psi^h\big) \mbox{~for all~} \psi^h\in
\mathcal{V}_{1,0}^h \label{10/4:ct1*}
\end{align}
has a unique solution $v^h\in \mathcal{V}_{1,\diamond}^h$ with $\gamma v^h = g$. Furthermore, the inequality
\begin{align}
\big\|v^h\big\|_{H^1(\Omega)}\le C_\mathcal{D} \left( \left\|f\right\|_{L^2(\Omega)}+
\left\|g\right\|_{H^{1/2}(\partial\Omega)} \right) \label{18/5:ct1*}
\end{align}
is satisfied. The map $\mathcal{D}^h: L^2(\Omega) \rightarrow
\mathcal{V}_{1,\diamond}^h$ from each $f \in L^2(\Omega)$
to the unique solution $v^h =: \mathcal{D}^h_f g$ of \eqref{10/4:ct1*}
is called the  {\it discrete Dirichlet operator}.
\end{proposition}
We now can introduce the strictly convex, discrete cost function
$$\Upsilon_{\rho,\delta}^h (f):=
\mathcal{J}_\delta^h(f) + \rho  \left \| f-f^*
\right \|^2_{L^2(\Omega)} \mbox{~with~} \mathcal{J}_\delta^h (f):= \int_\Omega Q\nabla \left(\mathcal{N}^h_fj_\delta-
\mathcal{D}^h_fg_\delta\right) \cdot \nabla \left(\mathcal{N}^h_fj_\delta-
\mathcal{D}^h_fg_\delta\right).$$
\begin{theorem}\label{Projection}
The problem
$$
\min_{f\in L^2(\Omega)} \Upsilon_{\rho,\delta}^h (f) \eqno \left(
\mathcal{P}_{\rho,\delta}^h \right)
$$
attains a unique minimizer $f$ which satisfies the equation
\begin{align}
f -f^* =-
\frac{1}{\rho} \left( \mathcal{N}^h_fj_\delta - \mathcal{D}^h_fg_\delta\right).\label{10-5-16ct3*}
\end{align}
\end{theorem}

\begin{remark}
Since $\mathcal{N}^h_fj_\delta$ and $\mathcal{D}^h_fg_\delta$ are both in $\mathcal{V}^1_h$, so is $f$, provided that $f^*\in \mathcal{V}^1_h$. Thus, taking this into account, a discretization of the set $L^2(\Omega)$ can be avoided.
\end{remark}

\begin{proof}[Proof of Theorem \ref{Projection}]
One can see easily that the problem $\left(\mathcal{P}_{\rho,\delta}^h \right)$ has a unique solution. It remains to show  \eqref{10-5-16ct3*}.
Let $f\in L^2(\Omega)$ be the minimizer to $\left(
\mathcal{P}_{\rho,\delta}^h \right)$. The first-order optimality condition yields that
$
{\Upsilon_{\rho,\delta}^h}'(f) \xi = {\mathcal{J}_\delta^h}' (f) \xi + 2\rho (\xi, f-f^*) = 0
$
for all $\xi \in L^2(\Omega)$. A short computation shows
$
{\mathcal{J}_\delta^h}'(f)(\xi) = 2\left(  \xi, \mathcal{N}^h_fj_\delta-
\mathcal{D}^h_fg_\delta \right)
$
and so obtain
$
\left(\xi, \frac{1}{\rho} \left(
\mathcal{N}^h_fj_\delta - \mathcal{D}^h_fg_\delta\right)  + f -f^*\right) = 0
$
for all $\xi \in L^2(\Omega)$,
which finishes the proof.
\end{proof}

From now on $C$ is a generic positive constant which is independent of the mesh size $h$ of $\mathcal{T}^h$, the noise level $\delta$ and the regularization parameter $\rho$. Before presenting the convergence of finite element approximations we here state some auxiliary results.

\begin{lemma}\label{moli.data}
A projection operator $\Pi^h_\diamond: L^1(\Omega) \rightarrow \mathcal{V}^h_{1,\diamond}$ exists such that
$$\Pi^h_\diamond\varphi^h = \varphi^h \mbox{~for all~} \varphi^h \in \mathcal{V}^h_{1,\diamond}  \mbox{~and~} \Pi^h_\diamond \big(H^1_0(\Omega)\big) \subset \mathcal{V}^h_{1,0} \subset \mathcal{V}^h_{1,\diamond}.$$
Furthermore, it satisfies the properties
\begin{equation}\label{23/10:ct2*}
\lim_{h\to 0} \big\| \vartheta - \Pi^h_\diamond \vartheta
\big\|_{H^1(\Omega)} =0 \enskip \mbox{~for all~} \vartheta \in H^1_\diamond(\Omega)
\end{equation}
and
\begin{equation}\label{23/5:ct1*}
\big\| \vartheta - \Pi^h_\diamond \vartheta \big\|_{H^1(\Omega)} \le
Ch \| \vartheta\|_{H^2(\Omega)} \mbox{~for all~} \vartheta \in H^1_\diamond(\Omega)\cap H^2(\Omega).
\end{equation}
\end{lemma}

\begin{proof}
Let $\Pi^h: L^1(\Omega) \rightarrow \mathcal{V}^h_1$ be the Clement's mollification interpolation operator, see \cite{Clement} and some generalizations \cite{Bernardi1,Bernardi2,scott_zhang}.
We then define the operator
$$\Pi^h_\diamond\vartheta := \Pi^h\vartheta -\frac{1}{|\partial\Omega|}\int_{\partial\Omega} \gamma \Pi^h\vartheta \in \mathcal{V}^h_{1,\diamond}, \quad \forall \vartheta\in L^1(\Omega)$$
which has the properties \eqref{23/10:ct2*} and \eqref{23/5:ct1*}. The proof is completed.
\end{proof}

On the basis of (\ref{23/10:ct2*}) and (\ref{23/5:ct1*}) we introduce for each $\Phi \in H^1_\diamond(\Omega)$
\begin{align}\label{3-6-16ct8}
\varrho^h_\Phi := \big\| \Phi - \Pi^h_\diamond \Phi \big\|_{H^1(\Omega)}.
\end{align}
We note that $\lim_{h \to 0} \varrho^h_\Phi = 0$ and
\begin{align}\label{21-5-16ct2*}
0 \le \varrho^h_\Phi \le Ch
\end{align}
in case $\Phi \in H^2(\Omega)$. Furthermore, let $(f,j,g)\in L^2(\Omega) \times H^{-1/2}(\partial\Omega) \times H^{1/2}_\diamond(\partial\Omega)$ be fixed, we denote by
\begin{align}\label{3-6-16ct9}
\alpha^h_{f,j}= \big\|\mathcal{N}^{h}_{f} j - \mathcal{N}_{f} j \big\|_{H^1(\Omega)} \mbox{~and~} \beta^h_{f,g}= \big\|\mathcal{D}^{h}_{f} g - \mathcal{D}_{f} g \big\|_{H^1(\Omega)}.
\end{align}
Then
$
\lim_{h\to 0} \alpha^h_{f,j} = \lim_{h\to 0}\beta^h_{f,g} =0.
$
In particular, if $\mathcal{N}_{f} j \in H^2(\Omega)$ and $\mathcal{D}_{f} g \in H^2(\Omega)$, the error estimates
\begin{align}\label{21-5-16ct2}
\alpha^h_{f,j} \le Ch
\mbox{~and~} \beta^h_{f,g} \le Ch
\end{align}
are satisfied (cf.\ \cite{Brenner_Scott,Ciarlet}).

\begin{lemma}\label{6-5-16ct5}
Let $(f_1,j_1,g_1)$ and $(f_2,j_2,g_2)$ be arbitrary in $ L^2(\Omega) \times H^{-1/2}(\partial\Omega) \times H^{1/2}_\diamond(\partial\Omega)$. Then the estimates
\begin{align}\label{7-5-16ct1}
\left\|\mathcal{N}^{h}_{f_1} j_1 - \mathcal{N}^{h}_{f_2} j_2 \right\|_{H^1(\Omega)} \le C_\mathcal{N} \left( \left\|f_1 - f_2\right\|_{L^2(\Omega)} + \left\|j_1 -j_2 \right\|_{H^{-1/2}(\partial\Omega)}\right)
\end{align}
and
\begin{align}\label{7-5-16ct2}
\left\|\mathcal{D}^{h}_{f_1} g_1 - \mathcal{D}^{h}_{f_2} g_2 \right\|_{H^1(\Omega)} \le C_\mathcal{D} \left( \left\|f_1 - f_2\right\|_{L^2(\Omega)} + \left\|g_1 -g_2 \right\|_{H^{1/2}(\partial\Omega)}\right)
\end{align}
hold for all $h \ge 0$.
\end{lemma}

\begin{proof}
According to the definition of the discrete Neumann operator, we have for all $\varphi^{h}\in
\mathcal{V}_{1,\diamond}^{h}$ that
$\int_\Omega  Q \nabla \mathcal{N}^{h}_{f_i} j_i \cdot \nabla \varphi^{h} =
\left\langle j_i,\gamma \varphi^{h}\right\rangle + \left( f_i,\varphi^{h}\right) \mbox{~with~} i=1,2$.
Thus, $\Phi^{h}_{\mathcal{N}} := \mathcal{N}^{h}_{f_1} j_1 - \mathcal{N}^{h}_{f_2} j_2$ is the unique solution to the variational problem
$\int_\Omega  Q \nabla \Phi^{h}_{\mathcal{N}} \cdot \nabla \varphi^{h} =
\left\langle j_1-j_2,\gamma \varphi^{h}\right\rangle + \left( f_1-f_2,\varphi^{h}\right) $
for all $\varphi^{h}\in
\mathcal{V}_{1,\diamond}^{h}$
and so that \eqref{7-5-16ct1} is satisfied. Likewise, we also obtain \eqref{7-5-16ct2}.
The proof is completed.
\end{proof}

\begin{lemma}\label{6-5-16ct6}
Let $\left(\mathcal{T}^{h_n}\right)$ be a sequence of triangulations with
$\limn h_n = 0$. Assume that $\left( j_{\delta_n}, g_{\delta_n}\right)$ is a sequence in $H^{-1/2}(\partial\Omega) \times H^{1/2}_\diamond(\partial\Omega)$ convergent to $\left( j_\delta, g_\delta\right) $ in the  $H^{-1/2}(\partial\Omega) \times H^{1/2}(\partial\Omega)$-norm and $(f_n)$ is a sequence in $L^2(\Omega)$ weakly convergent in $L^2(\Omega)$ to $f$, then there holds the inequality
\begin{align} \label{20-5-16ct3}
\liminfn \mathcal{J}^{h_n}_{\delta_n} (f_n) \ge \mathcal{J}_{\delta} (f).
\end{align}
\end{lemma}

\begin{proof}
The proof is based upon the mollification operator introduced in Lemma \ref{moli.data} together with standard arguments, therefore omitted here.
\end{proof}

We now show the convergence of finite element approximations to the identification problem.

\begin{theorem}\label{stability2}
Assume that
$\limn h_n = 0$ and $(\delta_n)$ and $(\rho_n)$ any positive sequences such that
\begin{align}\label{31-8-16ct1}
\rho_n \to 0, ~\frac{\delta_n}{\sqrt{\rho_n}} \to 0,~ \frac{\alpha^{h_n}_{f^\dag,j^\dag}}{\sqrt{\rho_n}} \to 0 \mbox{~and~} \frac{\beta^{h_n}_{f^\dag,g^\dag}}{\sqrt{\rho_n}} \to 0\mbox{~as~} n\to\infty,
\end{align}
where $\alpha^{h_n}_{f^\dag,j^\dag}$ and $\beta^{h_n}_{f^\dag,g^\dag}$ are defined by \eqref{3-6-16ct9}. Furthermore, assume that
$\left( j_{\delta_n}, g_{\delta_n}\right)$ is a sequence in $H^{-1/2}(\partial\Omega) \times H^{1/2}_\diamond(\partial\Omega)$ satisfying
$$\big\|j_{\delta_n} -j^\dag\big\|_{H^{-1/2}(\partial\Omega)} + \big\|g_{\delta_n} - g^\dag\big\|_{H^{1/2}(\partial\Omega)} \le \delta_n$$
and $f_n := f^{h_n}_{\rho_n,\delta_n}$ is the unique minimizer of $\left(\mathcal{P}^{h_n}_{\rho_n,\delta_n}\right)$ for each $n\in N$. Then:

(i) The sequence $(f_n)$ converges in the $L^2(\Omega)$-norm to $f^\dag$.

(ii) The corresponding state sequences $\left( \mathcal{N}^{h_n}_{f_n}j_{\delta_n}\right) $ and $\left( \mathcal{D}^{h_n}_{f_n}g_{\delta_n}\right) $ converge in the $H^1(\Omega)$-norm to the unique weak solution $\Phi^\dag = \Phi^\dag \big(f^\dag,j^\dag,g^\dag\big)$ of the boundary value problem \eqref{17-5-16ct1}--\eqref{17-5-16ct3}.
\end{theorem}

Before going to prove the theorem, we make the following short remark.

\begin{remark}\label{regu}
In case the weak solution $\Phi^\dag = \Phi^\dag \big(f^\dag,j^\dag,g^\dag\big)$ of \eqref{17-5-16ct1}--\eqref{17-5-16ct3} belonging to $H^2(\Omega)$, the estimate \eqref{21-5-16ct2} shows that  $0\le \alpha^{h_n}_{f^\dag,j^\dag},~\beta^{h_n}_{f^\dag,g^\dag} \le Ch_n$. Therefore, in view of \eqref{31-8-16ct1}, the above convergences (i) and (ii) are obtained if the sequence $(\rho_n)$ is chosen such that
\begin{align*}
\rho_n\rightarrow 0, ~\frac{\delta_n}{\sqrt{\rho_n}}
\rightarrow 0 \mbox{~and~} \frac{h_n}{\sqrt{\rho_n}}
\rightarrow 0 \mbox{~as~} n\to\infty.
\end{align*}
By regularity theory for elliptic boundary value problems, the regularity assumption $\Phi^\dag \in H^2(\Omega)$ is satisfied if the diffusion matrix $Q\in {C^{0,1}(\Omega)}^{d\times d}$, $j^\dag \in H^{1/2}(\partial\Omega)$, $g^\dag \in H^{3/2}(\partial\Omega)$ and either $\partial\Omega$ is smooth of the class $C^{0,1}$ or the domain $\Omega$ is convex (see, for example, \cite{Grisvad,Troianiello}).
\end{remark}

\begin{proof}[Proof of Theorem \ref{stability2}]
We have from the optimality of $f_n$ that
\begin{align}\label{6-5-16ct8}
\mathcal{J}^{h_n}_{\delta_n} \left(f_n\right) + \rho_n \|f_n-f^*\|^2_{L^2(\Omega)} &\le \mathcal{J}^{h_n}_{\delta_n} (f^\dag) + \rho_n \|f^\dag-f^*\|^2_{L^2(\Omega)}.
\end{align}
Since at $f^\dag$ there holds the equation $\mathcal{N}_{f^\dag} j^\dag = \mathcal{D}_{f^\dag} g^\dag$, we infer from Lemma \ref{6-5-16ct5} that
\begin{align}\label{6-5-16ct9}
\mathcal{J}^{h_n}_{\delta_n} (f^\dag)
\le C \left(\delta^2_n + \left(\alpha^{h_n}_{f^\dag, j^\dag}\right)^2 + \left(\beta^{h_n}_{f^\dag, g^\dag}\right)^2 \right)
\end{align}
which implies from \eqref{6-5-16ct8} that
$
\limn \mathcal{J}^{h_n}_{\delta_n} \left(f_n\right) =0
$
and, by the assumption \eqref{31-8-16ct1},
\begin{align}\label{6-5-16ct11}
\limsupn \|f_n-f^*\|^2_{L^2(\Omega)} \le  \big\|f^\dag-f^*\big\|^2_{L^2(\Omega)}.
\end{align}
So that the sequence $\left( f_n\right) $ is bounded in the $L^2(\Omega)$-norm. A subsequence not relabelled and an element $\widehat{f} \in L^2(\Omega)$ exist such that $\left( f_n\right)$ converges weakly in $L^2(\Omega)$ to $\widehat{f}$ and
\begin{align}\label{6-5-16ct12}
\big\|\widehat{f}-f^*\big\|^2_{L^2(\Omega)} \le \liminfn \|f_n-f^*\|^2_{L^2(\Omega)}.
\end{align}
For any $f\in L^2(\Omega)$ we denote by
$
\mathcal{J}_0(f) := \int_\Omega Q\nabla\left(\mathcal{N}_f j^\dag - \mathcal{D}_f g^\dag\right) \cdot \nabla\left(\mathcal{N}_f j^\dag - \mathcal{D}_f g^\dag\right).
$
By \eqref{coercivity}, we have
$
\left\|\mathcal{N}_{\widehat{f}} j^\dag - \mathcal{D}_{\widehat{f}} g^\dag\right\|^2_{H^1(\Omega)} \le \frac{1+ C_\Omega}{C_\Omega \underline{q}} \mathcal{J}_0 \big(\widehat{f}\big) \le \frac{1+ C_\Omega}{C_\Omega \underline{q}}\liminfn \mathcal{J}^{h_n}_{\delta_n}\left( f_n\right) = 0
$, here we used Lemma \ref{6-5-16ct6}.
Thus,
$\mathcal{N}_{\widehat{f}} j^\dag = \mathcal{D}_{\widehat{f}} g^\dag$ which infers $\widehat{f} \in \mathcal{I} \left(j^\dag,g^\dag\right).$
Now we show $\widehat{f} = f^\dag$ and the sequence $\left( f_n\right)$ converges to $\widehat{f}$ in the $L^2(\Omega)$-norm. By the definition of the $f^*$-minimum-norm solution and \eqref{6-5-16ct11}--\eqref{6-5-16ct12}, we get that
$$
\big\|f^\dag-f^*\big\|^2_{L^2(\Omega)} \le \big\|\widehat{f}-f^*\big\|^2_{L^2(\Omega)} \le \liminfn \|f_n-f^*\|^2_{L^2(\Omega)} \le \limsupn \|f_n-f^*\|^2_{L^2(\Omega)} \le \big\|f^\dag-f^*\big\|^2_{L^2(\Omega)}
$$
and so that
$
\big\|f^\dag-f^*\big\|^2_{L^2(\Omega)} = \big\|\widehat{f}-f^*\big\|^2_{L^2(\Omega)} = \limn \|f_n-f^*\|^2_{L^2(\Omega)}.
$
By the uniqueness of the minimum-norm solution and the sequence $\left( f_n\right)$ weakly converging in $L^2(\Omega)$ to $\widehat{f}$, we conclude that $\widehat{f} = f^\dag$ and the sequence $\left( f_n\right)$ in fact converges in the $L^2(\Omega)$-norm to $\widehat{f}$. Finally, we show the sequences $( \mathcal{N}^{h_n}_{f_n}j_{\delta_n}) $ and $( \mathcal{D}^{h_n}_{f_n}g_{\delta_n}) $ converge to $\Phi^\dag = \mathcal{N}_{f^\dag} j^\dag = \mathcal{D}_{f^\dag} g^\dag$ in the $H^1(\Omega)$-norm. Indeed,
by Lemma \ref{6-5-16ct5}, we obtain that
\begin{align*}
\left\|\mathcal{N}^{h_n}_{f_n}j_{\delta_n} - \mathcal{N}_{f^\dag} j^\dag \right\|_{H^1(\Omega)} &\le \left\|\mathcal{N}^{h_n}_{f_n}j_{\delta_n} - \mathcal{N}^{h_n}_{f^\dag} j^\dag \right\|_{H^1(\Omega)} + \left\|\mathcal{N}^{h_n}_{f^\dag} j^\dag - \mathcal{N}_{f^\dag} j^\dag\right\|_{H^1(\Omega)}\\
& \le C \left( \left\|j_{\delta_n} -j^\dag \right\|_{H^{-1/2}(\partial\Omega)} + \left\|f_n - f^\dag\right\|_{L^2(\Omega)} + \alpha^{h_n}_{f^\dag, j^\dag}\right) \rightarrow 0 \mbox{~as~} n \to \infty.
\end{align*}
Similarly, we also get
$
\left\|\mathcal{D}^{h_n}_{f_n}g_{\delta_n} - \mathcal{D}_{f^\dag} g^\dag \right\|_{H^1(\Omega)}
 \le C \left( \left\|g_{\delta_n} -g^\dag \right\|_{H^{1/2}(\partial\Omega)} + \left\|f_n - f^\dag\right\|_{L^2(\Omega)} + \beta^{h_n}_{f^\dag, g^\dag}\right) \rightarrow 0
$
as $n$ tends to $\infty$, which finishes the proof.
\end{proof}

\section{Convergence rates for the discretized problem}\label{tdht}

We are now in a position to state the main theorem on convergence rates for the general case of finite element discretized regularized solutions with noise level ($\delta>0$ and $h>0$). The source condition (\ref{error5}) will play a prominent role in this context.

\begin{theorem}\label{con.rate}
Assume that the condition (\ref{error5}) is fulfilled.
Then, we have
the error estimate and convergence rate
\begin{align} \label{29-6-15ct3}
\Big\| \mathcal{N}^h_{f^h} j_\delta &- \mathcal{D}^h_{f^h} g_\delta\Big\|^2_{H^1(\Omega)} + \rho\big\|f^h - f^\dag \big\|^2_{L^2(\Omega)} \notag\\
& =\mathcal{O} \left(\delta^2 + \left( \alpha^h_{f^\dag, j^\dag}\right)^2 + \left( \beta^h_{f^\dag, g^\dag}\right)^2 + \rho \varrho^h_{\mathcal{N}_w j^\dag} + \rho\varrho^h_{\mathcal{N}_{f^\dag} j^\dag} + \rho\varrho^h_{\mathcal{D}_0 \gamma \mathcal{N}_w j^\dag - g^\dag} +\delta\rho +\rho^2 \right),
\end{align}
where $f^h := f_{\rho,\delta}^h$ is the unique minimizer of $\left(
\mathcal{P}_{\rho,\delta}^h \right)$ and $\mathcal{D}_0 \gamma \mathcal{N}_w j^\dag - g^\dag$ is the unique weak solution to the Dirichlet problem
\begin{align*}
-\nabla \cdot (Q\nabla v) = 0 \mbox{~in~} \Omega \mbox{~and~} v = \gamma \mathcal{N}_w j^\dag - g^\dag \mbox{~on~} \partial\Omega
\end{align*}
and $\alpha^h_{f^\dag, j^\dag}$, $\beta^h_{f^\dag, g^\dag}$, $\varrho^h_{\mathcal{N}_w j^\dag}$, $\varrho^h_{\mathcal{N}_{f^\dag} j^\dag}$ and $\varrho^h_{\mathcal{D}_0 \gamma \mathcal{N}_w j^\dag - g^\dag}$ come from \eqref{3-6-16ct8} and \eqref{3-6-16ct9}.
\end{theorem}

\begin{remark}
In case (cf.\ Remark \ref{regu})
$
\mathcal{N}_{f^\dag} j^\dag ,\mathcal{N}_w j^\dag, \mathcal{D}_0 \gamma \mathcal{N}_w j^\dag - g^\dag\in H^2(\Omega),
$
by \eqref{21-5-16ct2*} and \eqref{21-5-16ct2}, we have
$$0\le \alpha^h_{f^\dag, j^\dag}, \beta^h_{f^\dag, g^\dag}, \varrho^h_{\mathcal{N}_w j^\dag}, \varrho^h_{\mathcal{N}_{f^\dag} j^\dag}, \varrho^h_{\mathcal{D}_0 \gamma \mathcal{N}_w j^\dag - g^\dag} \le Ch$$
and so that the following convergence rate is obtained
\begin{align*}
\left\| \mathcal{N}^h_{f^h} j_\delta - \mathcal{D}^h_{f^h} g_\delta\right\|^2_{H^1(\Omega)}  + \rho\big\|f^h - f^\dag \big\|^2_{L^2(\Omega)} =\mathcal{O} \left(\delta^2 + h^2 + h\rho + \delta\rho + \rho^2 \right).
\end{align*}
\end{remark}

\begin{remark}
Let $\Phi^\dag = \Phi^\dag \big(f^\dag,j^\dag,g^\dag\big)$ be the weak solution of \eqref{17-5-16ct1}--\eqref{17-5-16ct3}. Then the convergence rate
\begin{align*}
&\Big\|\mathcal{N}^h_{f^h}j_\delta - \Phi^\dag \Big\|^2_{H^1(\Omega)} + \left\|\mathcal{D}^h_{f^h}g_\delta - \Phi^\dag \right\|^2_{H^1(\Omega)}  \notag\\
& =\mathcal{O} \left(\delta^2\rho^{-1} + \left( \alpha^h_{f^\dag, j^\dag}\right)^2\rho^{-1} + \left( \beta^h_{f^\dag, g^\dag}\right)^2\rho^{-1} + \varrho^h_{\mathcal{N}_w j^\dag} + \varrho^h_{\mathcal{N}_{f^\dag} j^\dag} + \varrho^h_{\mathcal{D}_0 \gamma \mathcal{N}_w j^\dag - g^\dag} +\delta +\rho + \alpha^h_{f^\dag, j^\dag} + \beta^h_{f^\dag, g^\dag}\right)
\end{align*}
is also established. Indeed, the desired equation directly follows from  \eqref{29-6-15ct3} and the following inequalities
\begin{align*}
\left\|\mathcal{N}^{h}_{f^h}j_{\delta} - \mathcal{N}_{f^\dag} j^\dag \right\|_{H^1(\Omega)}
& \le C \left( \big\|j_{\delta} -j^\dag \big\|_{H^{-1/2}(\partial\Omega)} + \big\|f^h - f^\dag\big\|_{L^2(\Omega)} + \alpha^{h}_{f^\dag, j^\dag}\right) \\
& \le C \left( \delta + \big\|f^h - f^\dag\big\|_{L^2(\Omega)} + \alpha^{h}_{f^\dag, j^\dag}\right)
\end{align*}
and
$
\left\|\mathcal{D}^{h}_{f^h}g_{\delta} - \mathcal{D}_{f^\dag} g^\dag \right\|_{H^1(\Omega)}  \le C \left( \delta + \big\|f^h - f^\dag\big\|_{L^2(\Omega)} + \beta^{h}_{f^\dag, g^\dag}\right),
$
here we used Lemma \ref{6-5-16ct5}.
\end{remark}

\begin{proof}[Proof of Theorem \ref{con.rate}]
In view of \eqref{6-5-16ct9} we first have that
$\mathcal{J}^h_\delta \big(f^\dag\big) \le C\left(\delta^2 + \left( \alpha^h_{f^\dag, j^\dag}\right)^2 + \left( \beta^h_{f^\dag, g^\dag}\right)^2 \right).$
The optimality of $f^h$ yields
$
\mathcal{J}^h_\delta \big(f^h\big) + \rho\big\|f^h-f^*\big\|^2_{L^2(\Omega)} \le \mathcal{J}^h_\delta \big(f^\dag\big) +\rho \big\|f^\dag-f^*\big\|^2_{L^2(\Omega)}.
$
This gives
\begin{align}\label{2-6-16ct2}
\mathcal{J}^h_\delta \big(f^h\big) + \rho\big\|f^h - f^\dag \big\|^2_{L^2(\Omega)} &\le \mathcal{J}^h_\delta \big(f^\dag\big) +\rho \left( \big\|f^\dag-f^*\big\|^2_{L^2(\Omega)} - \big\|f^h-f^*\big\|^2_{L^2(\Omega)} + \big\|f^h - f^\dag \big\|^2_{L^2(\Omega)}\right) \notag\\
&\le C\left(\delta^2 + \left( \alpha^h_{f^\dag, j^\dag}\right)^2 + \left( \beta^h_{f^\dag, g^\dag}\right)^2 \right) + 2\rho\left( f^\dag-f^*, f^\dag-f^h\right).
\end{align}
Since $\mathcal{N}_{f^\dag} j^\dag = \mathcal{D}_{f^\dag} g^\dag$,
it follows that
\begin{align}\label{2-6-16ct3}
\left( f^\dag-f^*, f^\dag-f^h\right)
= \left( f^\dag-f^h, \mathcal{N}_wj^\dag - \mathcal{N}_{f^\dag} j^\dag \right) + \left( f^\dag-f^h, \mathcal{D}_{f^\dag} g^\dag - \mathcal{D}_wg^\dag \right).
\end{align}
From \eqref{ct9}, we infer
\begin{align*}
\left( f^\dag, \mathcal{N}_wj^\dag - \mathcal{N}_{f^\dag} j^\dag \right)
= \int_\Omega Q\nabla \mathcal{N}_{f^\dag} j^\dag \cdot \nabla \left( \mathcal{N}_w j^\dag - \mathcal{N}_{f^\dag} j^\dag\right) - \left\langle j^\dag, \gamma \left( \mathcal{N}_w j^\dag - \mathcal{N}_{f^\dag} j^\dag\right)\right\rangle,
\end{align*}
and
\begin{align*}
\left( f^h, \mathcal{N}_wj^\dag - \mathcal{N}_{f^\dag} j^\dag \right)
= \int_\Omega Q\nabla \mathcal{N}_{f^h} j^\dag \cdot \nabla \left( \mathcal{N}_w j^\dag - \mathcal{N}_{f^\dag} j^\dag\right) - \left\langle j^\dag, \gamma \left( \mathcal{N}_w j^\dag - \mathcal{N}_{f^\dag} j^\dag\right)\right\rangle.
\end{align*}
This in turn implies
\begin{align}\label{2-6-16ct4}
&\left( f^\dag-f^h, \mathcal{N}_wj^\dag - \mathcal{N}_{f^\dag} j^\dag \right)=
\int_\Omega Q\nabla \left( \mathcal{N}_{f^\dag} j^\dag - \mathcal{N}_{f^h} j^\dag\right) \cdot \nabla \left( \mathcal{N}_w j^\dag - \mathcal{N}_{f^\dag} j^\dag\right) \notag\\
&= \int_\Omega Q\nabla \left( \mathcal{N}_{f^\dag} j^\dag - \mathcal{D}_{f^h} g^\dag\right) \cdot \nabla \left( \mathcal{N}_w j^\dag - \mathcal{N}_{f^\dag} j^\dag\right)  + \int_\Omega Q\nabla \left( \mathcal{D}_{f^h} g^\dag - \mathcal{N}_{f^h} j^\dag\right) \cdot \nabla \left( \mathcal{N}_w j^\dag - \mathcal{N}_{f^\dag} j^\dag\right).
\end{align}
Since $\gamma \left( \mathcal{D}_{f^\dag} g^\dag - \mathcal{D}_wg^\dag\right) =0$, it follows from \eqref{ct9*} that
\begin{align}\label{2-6-16ct5}
\left( f^\dag-f^h, \mathcal{D}_{f^\dag} g^\dag - \mathcal{D}_wg^\dag \right)
= \int_\Omega Q \nabla \left( \mathcal{D}_{f^\dag} g^\dag - \mathcal{D}_{f^h} g^\dag\right)  \cdot \nabla \left( \mathcal{D}_{f^\dag} g^\dag - \mathcal{D}_wg^\dag \right)
\end{align}
holds. We thus infer from \eqref{2-6-16ct3}--\eqref{2-6-16ct5} the identity
\begin{align}\label{2-6-16ct6}
&\left( f^\dag-f^*, f^\dag-f^h\right)= \int_\Omega Q\nabla \left( \mathcal{D}_{f^h} g^\dag - \mathcal{N}_{f^h} j^\dag\right) \cdot \nabla \left( \mathcal{N}_w j^\dag - \mathcal{N}_{f^\dag} j^\dag\right) \notag\\
&~\quad + \int_\Omega Q\nabla \left( \mathcal{N}_{f^\dag} j^\dag - \mathcal{D}_{f^h} g^\dag\right) \cdot \nabla \left( \mathcal{N}_w j^\dag - \mathcal{N}_{f^\dag} j^\dag\right) + \int_\Omega Q \nabla \left( \mathcal{D}_{f^\dag} g^\dag - \mathcal{D}_{f^h} g^\dag\right)  \cdot \nabla \left( \mathcal{D}_{f^\dag} g^\dag - \mathcal{D}_wg^\dag \right).
\end{align}
We note again that $\mathcal{N}_{f^\dag} j^\dag = \mathcal{D}_{f^\dag} g^\dag$ and $\gamma \left( \mathcal{D}_{f^\dag} g^\dag - \mathcal{D}_{f^h}g^\dag\right) =0$. Then, together with \eqref{ct9} and \eqref{ct9*}, the last two terms on the right hand side of \eqref{2-6-16ct6} satisfy
\begin{align*}
\int_\Omega Q\nabla \left( \mathcal{N}_{f^\dag} j^\dag - \mathcal{D}_{f^h} g^\dag\right) \cdot \nabla \left( \mathcal{N}_w j^\dag - \mathcal{N}_{f^\dag} j^\dag\right) + \int_\Omega Q \nabla \left( \mathcal{D}_{f^\dag} g^\dag - \mathcal{D}_{f^h} g^\dag\right)  \cdot \nabla \left( \mathcal{D}_{f^\dag} g^\dag - \mathcal{D}_wg^\dag \right)
=0.
\end{align*}
Thus, we obtain from \eqref{2-6-16ct6}
\begin{align*}
\left( f^\dag-f^*, f^\dag-f^h\right) = \int_\Omega Q\nabla \left( \mathcal{D}_{f^h} g^\dag - \mathcal{N}_{f^h} j^\dag\right) \cdot \nabla \left( \mathcal{N}_w j^\dag - \mathcal{N}_{f^\dag} j^\dag\right).
\end{align*}
Next, we abbreviate
$
W=\mathcal{N}_w j^\dag - \mathcal{N}_{f^\dag} j^\dag
$
and note
\begin{align}\label{2-6-16ct7}
\gamma W = \gamma \mathcal{N}_w j^\dag - g^\dag.
\end{align}
Then we get
\begin{align}\label{3-6-16ct1}
\left( f^\dag-f^*, f^\dag-f^h\right)
&= \int_\Omega Q\nabla \left( \mathcal{D}_{f^h} g^\dag - \mathcal{D}^h_{f^h} g^\dag\right) \cdot \nabla W - \int_\Omega Q\nabla \left( \mathcal{N}_{f^h} j^\dag - \mathcal{N}^h_{f^h} j^\dag\right) \cdot \nabla W\notag\\
&~\quad + \int_\Omega Q\nabla \left( \mathcal{D}^h_{f^h} g^\dag - \mathcal{D}^h_{f^h} g_\delta\right) \cdot \nabla W - \int_\Omega Q\nabla \left( \mathcal{N}^h_{f^h} j^\dag - \mathcal{N}^h_{f^h} j_\delta\right) \cdot \nabla W\notag\\
&~\quad + \int_\Omega Q\nabla \left( \mathcal{D}^h_{f^h} g_\delta - \mathcal{N}^h_{f^h} j_\delta\right) \cdot \nabla W := I_1+I_2+I_3.
\end{align}
To prepare the estimation of those three addends we start with writing
\begin{align*}
&\int_\Omega Q\nabla \left( \mathcal{D}_{f^h} g^\dag - \mathcal{D}^h_{f^h} g^\dag\right) \cdot \nabla W \\
&= \int_\Omega Q\nabla \left( \mathcal{D}_{f^h} g^\dag - \mathcal{D}^h_{f^h} g^\dag\right) \cdot \nabla \mathcal{D}_0 \gamma W + \int_\Omega Q\nabla \left( \mathcal{D}_{f^h} g^\dag - \mathcal{D}^h_{f^h} g^\dag\right) \cdot \nabla \Pi^h_\diamond \left( W- \mathcal{D}_0 \gamma W\right) \\
&~\quad + \int_\Omega Q\nabla \left( \mathcal{D}_{f^h} g^\dag - \mathcal{D}^h_{f^h} g^\dag\right) \cdot \nabla \left( W- \mathcal{D}_0 \gamma W - \Pi^h_\diamond \left( W- \mathcal{D}_0 \gamma W\right) \right).
\end{align*}
Since $\mathcal{D}_{f^h} g^\dag - \mathcal{D}^h_{f^h} g^\dag \in H^1_0(\Omega)$, we then get
\begin{align*}
\int_\Omega Q\nabla \left( \mathcal{D}_{f^h} g^\dag - \mathcal{D}^h_{f^h} g^\dag\right) \cdot \nabla \mathcal{D}_0 \gamma W &= \int_\Omega Q\nabla \mathcal{D}_0 \gamma W \cdot \nabla \left( \mathcal{D}_{f^h} g^\dag - \mathcal{D}^h_{f^h} g^\dag\right)= 0.
\end{align*}
Since
$\gamma \left( W- \mathcal{D}_0 \gamma W\right) = \gamma W - \gamma \mathcal{D}_0 \gamma W = \gamma W - \gamma W =0$,
we infer $\Pi^h_\diamond \left( W- \mathcal{D}_0 \gamma W\right) \in \mathcal{V}_{1,0}^h = \mathcal{V}_1^h \cap H^1_0(\Omega)$ and then obtain from \eqref{ct9*} and \eqref{10/4:ct1*} that
$
\int_\Omega Q\nabla \left( \mathcal{D}_{f^h} g^\dag - \mathcal{D}^h_{f^h} g^\dag\right) \cdot \nabla \Pi^h_\diamond \left( W- \mathcal{D}_0 \gamma W\right)=0
$
holds. Hence we have
\begin{align}\label{3-6-16ct2}
&\bigg| \int_\Omega Q\nabla \left( \mathcal{D}_{f^h} g^\dag - \mathcal{D}^h_{f^h} g^\dag\right) \cdot \nabla W \bigg|\notag\\
&= \left| \int_\Omega Q\nabla \left( \mathcal{D}_{f^h} g^\dag - \mathcal{D}^h_{f^h} g^\dag\right) \cdot \nabla \left( W- \mathcal{D}_0 \gamma W - \Pi^h_\diamond \left( W- \mathcal{D}_0 \gamma W\right) \right)\right| \notag\\
&\le C \left( \big\|f^h\big\|_{L^2(\Omega)} +  \big\|g^\dag\big\|_{H^{1/2}(\partial\Omega)}\right) \left( \varrho^h_{\mathcal{N}_w j^\dag} + \varrho^h_{\mathcal{N}_{f^\dag} j^\dag} + \varrho^h_{\mathcal{D}_0 \gamma \mathcal{N}_w j^\dag - g^\dag}\right),
\end{align}
where we use  \eqref{2-6-16ct7}. Similarly, since $\Pi^h_\diamond W \in \mathcal{V}^h_\diamond$ and by  \eqref{ct9} and \eqref{10/4:ct1}, we get
\begin{align}\label{3-6-16ct3}
\left|  \int_\Omega Q\nabla \left( \mathcal{N}_{f^h} j^\dag - \mathcal{N}^h_{f^h} j^\dag\right) \cdot \nabla W \right|
\le C \left( \big\|f^h\big\|_{L^2(\Omega)} +  \big\|j^\dag\big\|_{H^{-1/2}(\partial\Omega)}\right)\left( \varrho^h_{\mathcal{N}_w j^\dag} + \varrho^h_{\mathcal{N}_{f^\dag} j^\dag}\right).
\end{align}
Now we are in the position to estimate $I_1-I_3$. Combining \eqref{3-6-16ct2} with \eqref{3-6-16ct3}, we obtain
\begin{align}\label{3-6-16ct4}
\rho |I_1|
\le C\left(\delta^2 + \left( \alpha^h_{f^\dag, j^\dag}\right)^2 + \left( \beta^h_{f^\dag, g^\dag}\right)^2 + \rho \varrho^h_{\mathcal{N}_w j^\dag} + \rho\varrho^h_{\mathcal{N}_{f^\dag} j^\dag} + \rho\varrho^h_{\mathcal{D}_0 \gamma \mathcal{N}_w j^\dag - g^\dag} \right).
\end{align}
Now, using Lemma \ref{6-5-16ct5}, we arrive at
\begin{align}\label{3-6-16ct5}
\rho|I_2|
\le C \rho\left( \big\|g_\delta-g^\dag\big\|_{H^{1/2}(\partial\Omega)} + \big\|j_\delta-j^\dag\big\|_{H^{-1/2}(\partial\Omega)}\right) \le C\delta\rho.
\end{align}
Since for a.e. in $\Omega$ the matrix $Q(x)$ is positive definite, the root $Q(x)^{1/2}$ is then well defined. Thus, using the Cauchy-Schwarz inequality and Young's inequality, we estimate $I_3$ as
\begin{align}\label{3-6-16ct6}
\rho|I_3|
\le C \rho\left( \mathcal{J}^h_\delta \big(f^h\big)\right)^{1/2} \le C^2\rho^2 + \frac{1}{4} \mathcal{J}^h_\delta \big(f^h\big) \le C\rho^2 + \frac{1}{4} \mathcal{J}^h_\delta \big(f^h\big).
\end{align}
It follows from \eqref{3-6-16ct1} and \eqref{3-6-16ct4}--\eqref{3-6-16ct6} that
\begin{align*}
&2\rho \left( f^\dag-f^*, f^\dag-f^h\right) \\
&\le C\left( \delta^2 + \left( \alpha^h_{f^\dag, j^\dag}\right)^2 + \left( \beta^h_{f^\dag, g^\dag}\right)^2  + \rho \varrho^h_{\mathcal{N}_w j^\dag} + \rho\varrho^h_{\mathcal{N}_{f^\dag} j^\dag} + \rho\varrho^h_{\mathcal{D}_0 \gamma \mathcal{N}_w j^\dag - g^\dag} + \rho\delta + \rho^2\right) + \frac{1}{2} \mathcal{J}^h_\delta \big(f^h\big)
\end{align*}
holds, which together with \eqref{2-6-16ct2} implies
\begin{align}\label{3-6-16ct7}
&\frac{1}{2}\mathcal{J}^h_\delta \big(f^h\big) + \rho\big\|f^h - f^\dag \big\|^2_{L^2(\Omega)} \notag\\
&\le C\left(\delta^2 + \left( \alpha^h_{f^\dag, j^\dag}\right)^2 + \left( \beta^h_{f^\dag, g^\dag}\right)^2 + \rho \varrho^h_{\mathcal{N}_w j^\dag} + \rho\varrho^h_{\mathcal{N}_{f^\dag} j^\dag} + \rho\varrho^h_{\mathcal{D}_0 \gamma \mathcal{N}_w j^\dag - g^\dag} +\rho\delta +\rho^2 \right).
\end{align}
Since $\left\| \mathcal{D}^h_{f^h} g_\delta - \mathcal{N}^h_{f^h} j_\delta\right\|^2_{H^1(\Omega)} \le C \mathcal{J}^h_\delta \big(f^h\big),$
\eqref{29-6-15ct3} now directly follows from \eqref{3-6-16ct7}, which finishes the proof.
\end{proof}

\section{Conjugate gradient method and numerical test}\label{Numer_text}

In this section we will utilize the conjugate gradient (CG) method (see, for example, \cite{hager,kelley}) to find the minimizes of the strictly convex, discrete regularized problem $\left( \mathcal{P}_{\rho,\delta}^h \right)$.
Let
$
\nabla \Upsilon_{\rho,\delta}^h (f) = 2 \left(\mathcal{N}^h_fj_\delta-
\mathcal{D}^h_fg_\delta\right)
 +2\rho (f-f^*)
$
be the $L^2$-gradient of the cost function $\Upsilon_{\rho,\delta}^h$ at $f$ (see Proof of Theorem \ref{Projection}), where $f^*\in \mathcal{V}^h_1$. Then the sequence of iterates via this algorithm is generated by $f^0 \in L^2(\Omega) \cap\mathcal{V}^h_1$ and
$
f^{k+1} :=f^k + t^k d^k
$
for $k\ge 0$, where
\begin{align*}
d^k :=
\begin{cases}
-\nabla \Upsilon_{\rho,\delta}^h (f^k) &\mbox{~if~} k=0,\\
-\nabla \Upsilon_{\rho,\delta}^h (f^k) +\beta^kd^{k-1} &\mbox{~if~} k>0
\end{cases} \mbox{~with~} \beta^k := \frac{\| \nabla \Upsilon_{\rho,\delta}^h (f^k) \|^2}{\| \nabla \Upsilon_{\rho,\delta}^h (f^{k-1}) \|^2} \mbox{~and~}
t^k := \arg \min_{t \ge 0} \Upsilon_{\rho,\delta}^h (f^k + t d^k).
\end{align*}
A short computation shows that
\begin{align*}
t^k &= - \frac{\int_\Omega Q\nabla \left( \mathcal{N}^h_{d^k} 0 - \mathcal{D}^h_{d^k} 0 \right) \cdot \nabla \left( \mathcal{N}^h_{f^k} j_\delta - \mathcal{D}^h_{f^k} g_\delta\right) + \rho\left( d^k, f^k-f^*\right) }{\int_\Omega Q\nabla \left( \mathcal{N}^h_{d^k} 0 - \mathcal{D}^h_{d^k} 0 \right) \cdot \nabla \left( \mathcal{N}^h_{d^k} 0 - \mathcal{D}^h_{d^k} 0 \right) + \rho\left\|d^k\right\|^2_{L^2(\Omega)}}= - \frac{1}{2}\frac{\left( d^k, \nabla \Upsilon_{\rho,\delta}^h (f^k)\right)}{\left( d^k, \mathcal{N}^h_{d^k} 0 - \mathcal{D}^h_{d^k} 0 \right) + \rho\left\|d^k\right\|^2_{L^2(\Omega)}}.
\end{align*}
Consequently, the CG method then reads as follows:
giving an initial approximation $f^0\in \mathcal{V}^h_1$, number of iterations $N$ and a positive constants $\tau_1,\tau_2$. Computing
\begin{align*}
& \nabla \Upsilon_{\rho,\delta}^h (f^0) = 2 \left(\mathcal{N}^h_{f^0}j_\delta-
\mathcal{D}^h_{f^0}g_\delta\right) +2\rho (f^0-f^*),~  d^0 = -\nabla \Upsilon_{\rho,\delta}^h (f^0),~ t^0 = \frac{1}{2}\frac{\left\|d^0\right\|^2_{L^2(\Omega)}}{\left( d^0, \mathcal{N}^h_{d^0} 0 - \mathcal{D}^h_{d^0} 0 \right) + \rho\left\|d^0\right\|^2_{L^2(\Omega)}}
\end{align*}
and setting
\begin{align*}
f^1 = f^0 + t^0 d^0  \mbox{~and~} k=1, \quad \mbox{Tolerance}:= \big\|\nabla \Upsilon^h_{\rho,\delta} \big( f^k \big) \big\|_{L^2(\Omega)}-\tau_1-\tau_2\big\|\nabla \Upsilon^h_{\rho,\delta} \big( f^0 \big) \big\|_{L^2(\Omega)}.
\end{align*}

\begin{algorithm}[H]
\While{ $(\mbox{Tolerance}>0) ~\&~ (k\le N)$}
{
\begin{align*}
&\overline{r} = \left\| \nabla \Upsilon_{\rho,\delta}^h (f^{k-1}) \right\|^2_{L^2(\Omega)}, \quad r = \left\| \nabla \Upsilon_{\rho,\delta}^h (f^k) \right\|^2_{L^2(\Omega)}, \quad \beta^k = \dfrac{r}{\overline{r}},\\
& d^k = -\nabla \Upsilon_{\rho,\delta}^h (f^k) + \beta^k d^{k-1}, \quad t^k = -\frac{1}{2}\frac{\left( d^k, \nabla \Upsilon_{\rho,\delta}^h (f^k)\right)}{\left( d^k, \mathcal{N}^h_{d^k} 0 - \mathcal{D}^h_{d^k} 0 \right) + \rho\left\|d^k\right\|^2_{L^2(\Omega)}},\\
&f^{k+1} = f^k + t^k d^k,\\
&k:= k+1, \quad \mbox{Tolerance}:= \big\|\nabla \Upsilon^h_{\rho,\delta} \big( f^k \big) \big\|_{L^2(\Omega)}-\tau_1-\tau_2\big\|\nabla \Upsilon^h_{\rho,\delta} \big( f^0 \big) \big\|_{L^2(\Omega)}.
\end{align*}
}
\caption{CG iteration} \label{alg1}
\end{algorithm}

Below we illustrate the theoretical result with numerical examples. For this purpose we consider the the boundary value problem
\begin{align}
-\nabla \cdot \big(Q \nabla \Phi \big) &= f^\dag \mbox{~in~} \Omega := (-1,1) \times (-1,1),  \label{17-5-16ct1*}\\
Q \nabla \Phi \cdot \vec{n} &= j^\dag \mbox{~on~} \partial\Omega \mbox{~and~}
\Phi = g^\dag \mbox{~on~} \partial\Omega. \label{17-5-16ct3*}
\end{align}
We assume that entries of the known symmetric diffusion matrix $Q$ are discontinuous which are defined as
$q_{11}= 3\chi_{\Omega_{11}} + \chi_{\Omega\setminus\Omega_{11}}, ~
q_{12}= \chi_{\Omega_{12}}, ~
q_{22}= 4\chi_{\Omega_{22}} + 2\chi_{\Omega\setminus\Omega_{22}}$,
where $\chi_D$ is the characteristic function of the Lebesgue measurable set $D$ and
\begin{align*}
&\Omega_{11} := \left\{ (x_1, x_2) \in \Omega ~\big|~ |x_1| \le 1/2 \mbox{~and~} |x_2| \le 1/2 \right\},~ \Omega_{12} := \left\{ (x_1, x_2) \in \Omega ~\big|~ |x_1| + |x_2| \le 1/2 \right\} \mbox{~and~}\\
&\Omega_{22} := \left\{ (x_1, x_2) \in \Omega ~\big|~ x_1^2 + x_2^2 \le  1/4 \right\}.
\end{align*}

The identified source function $f^\dag \in L^2(\Omega)$ in \eqref{17-5-16ct1*} is assumed to be discontinuous and defined as
$$f^\dag = 2 \chi_{\Omega_{1}} - \chi_{\Omega_{2}} +\frac{5\pi}{7\pi-192} \chi_{\Omega \setminus (\Omega_1 \cup \Omega_2)},$$
where
\begin{align*}
&\Omega_1 := \left\{ (x_1, x_2) \in \Omega ~\big|~ 9(x_1 +1/2)^2 + 16(x_2-1/2)^2 \le 1\right\} \mbox{~and~}\\
&\Omega_2 := \left\{ (x_1, x_2) \in \Omega ~\big|~ (x_1-1/2)^2 + (x_2+1/2)^2 \le 1/16\right\}.
\end{align*}

For the discretization we divide the interval $(-1,1)$ into $\ell$ equal segments and so that the domain $\Omega = (-1,1)^2$ is divided into $2\ell^2$ triangles, where the diameter of each triangle is $h_{\ell} = \frac{\sqrt{8}}{\ell}$. In the minimization problem $\left(\mathcal{P}_{\rho,\delta}^{h} \right)$ we take $h=h_\ell$ and $\rho = \rho_\ell = 0.01h_\ell$. We use Algorithm \ref{alg1} which is described above for computing the numerical solution of the problem $\left(\mathcal{P}_{\rho_\ell,\delta_\ell}^{h_\ell} \right)$.
As an a-priori estimate and the initial approximation we choose $f^* := 0$ and $f^0(x) := \chi_{(0,1]\times[-1,1]} -\chi_{[-1,0]\times[-1,1]}$.

\begin{example}\label{30-6-17ct1}
In this first example $j^\dag \in H^{-1/2}(\partial\Omega)$ is chosen to be the piecewise constant function defined by
\begin{equation}\label{20-6-17ct5}
\begin{aligned}
j^\dag
&= \chi_{(0,1]\times\{-1\}} - \chi_{[-1,0]\times\{1\}} + 2\chi_{(0,1]\times\{1\}} -2\chi_{[-1,0]\times\{-1\}}\\
&~\quad +3\chi_{\{-1\}\times(-1,0]} - 3\chi_{\{1\}\times(0,1)} + 4\chi_{\{1\}\times(-1,0]}  - 4\chi_{\{-1\}\times(0,1)} .
\end{aligned}
\end{equation}
Then $g^\dag \in H^{1/2}_\diamond(\partial\Omega)$ is defined as
$g^\dag = \gamma \mathcal{N}_{f^\dag} j^\dag$. We mention that, to avoid a so-called inverse crime, we generate the given data on a finer grid than those used in the computations. For this purpose we first solve the problem \eqref{17-5-16ct1*} supplemented with the Neumann boundary condition in \eqref{17-5-16ct3*} on the very fine grid with $\ell = 128$, and then use this numerical approximation as substitute for $(j^\dag,g^\dag)$  in our computational considerations below.

For observations with noise we assume that
\begin{align}\label{3-7-17ct1}
\left( j_{\delta_{\ell}}, g_{\delta_{\ell}} \right) = \left( j^\dag+ \theta_\ell\cdot R_{j^\dag}, g^\dag+ \theta_\ell\cdot R_{g^\dag}\right) \quad \mbox{for some} \quad \theta_\ell>0 \quad \mbox{depending on} \quad \ell,
\end{align}
where $R_{j^\dag}$ and $R_{g^\dag}$ are $\partial M^{h_\ell}\times 1$-matrices of random numbers on the interval $(-1,1)$ which are generated by the MATLAB function ``rand'', and $\partial M^{h_\ell}$ is the number of boundary nodes of the triangulation $\mathcal{T}^{h_\ell}$. The measurement error is then computed as $\delta_\ell = \big\|j_{\delta_\ell} -j^\dag\big\|_{L^2(\partial\Omega)} + \big\|g_{\delta_\ell} -g^\dag\big\|_{L^2(\partial\Omega)}.$
To satisfy the condition $\delta_\ell\cdot\rho^{-1/2}_\ell\to 0$ as $\ell\to\infty$ in Theorem \ref{stability2} we below take $\theta_\ell = h_\ell \sqrt{\rho_\ell}$.

We start with the coarsest level $\ell =4$.
At each iteration $k$ we compute
$$\mbox{Tolerance}:= \big\|\nabla \Upsilon_{\rho_{\ell},\delta_{\ell}}^{h_{\ell}} \big( f^k_\ell \big) \big\|_{L^2(\Omega)}-\tau_1-\tau_2\big\|\nabla \Upsilon_{\rho_{\ell},\delta_{\ell}}^{h_{\ell}} \big( f^0_\ell \big) \big\|_{L^2(\Omega)},$$
where $\tau_1 := 10^{-6}h_\ell^{1/2}$ and $\tau_2 := 10^{-4}h_\ell^{1/2}$. Then the iteration was stopped if
$\mbox{Tolerance} \le 0$
or the number of iterations reached the maximum iteration count of 600. After obtaining the numerical solution of the first iteration process with respect to the coarsest level $\ell=4$, we use its interpolation on the next finer mesh $\ell=8$ as an initial approximation $f^0$ for the algorithm on this finer mesh, and proceed similarly on the preceding refinement levels.

Let $f_\ell$ be the function which is obtained at {\it the final iterate} of Algorithm \ref{alg1} corresponding to the refinement level $\ell$. Furthermore, let $\mathcal{N}^{h_\ell}_{f_\ell} j_{\delta_\ell}$ and $\mathcal{D}^{h_\ell}_{f_\ell} g_{\delta_\ell}$ denote {\it the computed numerical solution} to the Neumann and Dirichlet problem
\begin{align*}
-\nabla \cdot (Q\nabla u) = f_\ell \mbox{~in~} \Omega \mbox{~and~}
Q\nabla u\cdot \vec{n} = j_{\delta_\ell} \mbox{~on~} \partial\Omega
\quad\mbox{and}\quad
-\nabla \cdot (Q\nabla v) = f_\ell \mbox{~in~} \Omega \mbox{~and~}
v = g_\ell \mbox{~on~} \partial\Omega,
\end{align*}
respectively. The notations $\mathcal{N}^{h_\ell}_{ f^\dag}  j^\dag$ and $\mathcal{D}^{h_\ell}_{f^\dag} g^\dag$ of {\it the exact numerical solutions} are to be understood similarly. We use the following abbreviations for the errors
\begin{align*}
& L^2_f = \big\|f_\ell -  f^\dag\big\|_{L^2(\Omega)}, L^2_\mathcal{N} = \big\|\mathcal{N}^{h_\ell}_{f_\ell} j_{\delta_\ell} - \mathcal{N}^{h_\ell}_{ f^\dag}  j^\dag\big\|_{L^2(\Omega)},~ H^1_\mathcal{N} = \big\|\mathcal{N}^{h_\ell}_{f_\ell} j_{\delta_\ell} - \mathcal{N}^{h_\ell}_{ f^\dag}  j^\dag\big\|_{H^1(\Omega)} \mbox{~and~}\\
& L^2_\mathcal{D} = \big\|\mathcal{D}^{h_\ell}_{f_\ell} g_{\delta_\ell} - \mathcal{D}^{h_\ell}_{ f^\dag}  g^\dag\big\|_{L^2(\Omega)},~ H^1_\mathcal{D} = \big\|\mathcal{D}^{h_\ell}_{f_\ell} g_{\delta_\ell} - \mathcal{D}^{h_\ell}_{ f^\dag}  g^\dag\big\|_{H^1(\Omega)}.
\end{align*}
The numerical results are summarized in Table \ref{b1} and Table \ref{b2}, where we present the refinement level $\ell$, mesh size $h_\ell$ of the triangulation, regularization parameter $\rho_\ell$, measured noise $\delta_\ell$, number of iterations, value of tolerances and errors $L^2_f$, $L^2_{\mathcal{N}}$, $L^2_{\mathcal{D}}$, $H^1_{\mathcal{N}}$ and $H^1_{\mathcal{D}}$. Their experimental order of convergence (EOC) is presented in Table \ref{b3}, where
$\mbox{EOC}_\Theta := \dfrac{\ln \Theta(h_1) - \ln \Theta(h_2)}{\ln h_1 - \ln h_2}$
and $\Theta(h)$ is an error function of the mesh size $h$.

All figures presented correspond to $\ell = 64$. Figure \ref{h1} from left to right shows the computed numerical solution $f_\ell$ of the algorithm at the final $579^{\mbox{th}}$-iteration, and the differences $\mathcal{N}^{h_\ell}_{ f^\dag}  j^\dag - \mathcal{N}^{h_\ell}_{f_\ell} j_{\delta_\ell}$, ~$ \mathcal{D}^{h_\ell}_{f^\dag} g^\dag-\mathcal{D}^{h_\ell}_{f_\ell} g_{\delta_\ell}$ and $\mathcal{D}^{h_\ell}_{f_\ell} g_{\delta_\ell} - \mathcal{N}^{h_\ell}_{f_\ell} j_{\delta_\ell}$.
\begin{table}[H]
\begin{center}
\begin{tabular}{|c|l|l|l|l|l|}
\hline \multicolumn{6}{|c|}{ {\bf Convergence history} }\\
\hline
$\ell$ &\scriptsize $h_\ell$ &\scriptsize $\rho_\ell$ &\scriptsize $\delta_\ell$ &\scriptsize {\bf Iterate} &\scriptsize {\bf Tolerance} \\
\hline
4   &0.7071 & 0.7071e-2& 0.1916& 312&  -3.0822e-5\\
\hline
8  &0.3536 & 0.3536e-2& 9.3172e-2& 387& -1.2739e-6\\
\hline
16  &0.1766 & 0.1766e-2& 4.1174e-2& 461&   -1.4029e-6\\
\hline
32  &8.8388e-2 & 0.8839e-3& 2.0932e-2& 505&  -1.8559e-7\\
\hline
64  &4.4194e-2 & 0.4419e-3& 7.2765e-3& 579& -7.3540e-9\\
\hline
\end{tabular}
\caption{Refinement level $\ell$, mesh size $h_\ell$ of the triangulation, regularization parameter $\rho_\ell$,  measurement noise $\delta_\ell$, number of iterates and value of Tolerance.}
\label{b1}
\end{center}
\end{table}
\begin{table}[H]
\begin{center}
\begin{tabular}{|c|l|l|l|l|l|}
 \hline \multicolumn{6}{|c|}{ {\bf Convergence history} }\\
 \hline
$\ell$ &\scriptsize $L^2_f$ &\scriptsize  $L^2_{\mathcal{N}}$ &\scriptsize $L^2_{\mathcal{D}}$ &\scriptsize $H^1_{\mathcal{N}}$ &\scriptsize $H^1_{\mathcal{D}}$ \\
\hline
4&    0.5215    & 2.0441e-2  & 2.0396e-2  & 6.9952e-2  & 6.9713e-2 \\
\hline
8&   0.3309    & 6.3175e-3  & 6.3083e-3  & 3.1374e-2  & 3.1311e-2 \\
\hline
16&   0.1915 & 2.0132e-3  & 2.0122e-3  & 1.7276e-2  & 1.7243e-2 \\
\hline
32&   0.1073 & 5.5434e-4  & 5.5426e-4  & 8.9136e-3  & 8.9130e-3 \\
\hline
64&   5.2568e-2 &  1.4669e-4 & 1.4666e-4  & 3.9352e-3  & 3.9347e-3 \\
\hline
\end{tabular}
\caption{Errors  $L^2_f$, $L^2_{\mathcal{N}}$, $L^2_{\mathcal{D}}$, $H^1_{\mathcal{N}}$ and $H^1_{\mathcal{D}}$.}
\label{b2}
\end{center}
\end{table}
\begin{table}[H]
\begin{center}
\begin{tabular}{|c|l|l|l|l|l|}
 \hline \multicolumn{6}{|c|}{ {\bf Experimental order of convergence} }\\
 \hline
$\ell$ &\scriptsize {\bf EOC$_{L^2_f}$} &\scriptsize {\bf EOC$_{L^2_{\mathcal{N}}}$} &\scriptsize {\bf EOC$_{L^2_{\mathcal{D}}}$} &\scriptsize {\bf EOC$_{H^1_{\mathcal{N}}}$} &\scriptsize {\bf EOC$_{H^1_{\mathcal{D}}}$}\\
\hline
4 & -- & -- &-- &-- &--\\
\hline
8  & 0.6563 & 1.6940 & 1.6930 &  1.1568 & 1.1548\\
\hline
16  & 0.7891 & 1.6499 & 1.6485 &  0.8608 & 0.8607\\
\hline
32  & 0.8357 & 1.8606 & 1.8601 &  0.9547 & 0.9520\\
\hline
64  & 1.0294 & 1.9180 & 1.9181 &  1.1796 & 1.1797\\
\hline
{\bf Mean of EOC}
& 0.8276 & 1.7806& 1.7799 & 1.0380  & 1.0368\\\hline
\end{tabular}
\caption{Experimental order of convergence between finest and coarsest
level for $L^2_f$, $L^2_{\mathcal{N}}$, $L^2_{\mathcal{D}}$, $H^1_{\mathcal{N}}$ and $H^1_{\mathcal{D}}$.}
\label{b3}
\end{center}
\end{table}
\begin{figure}[H]
\begin{center}
\includegraphics[scale=0.2]{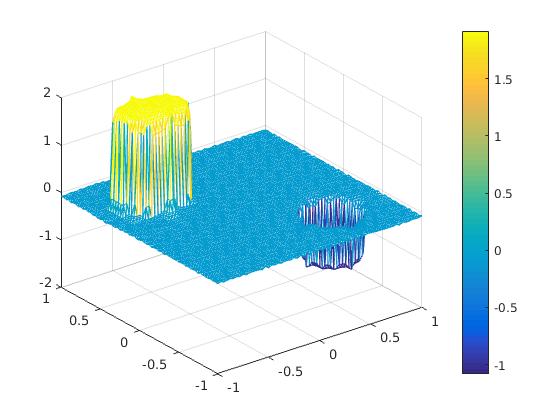}
\includegraphics[scale=0.2]{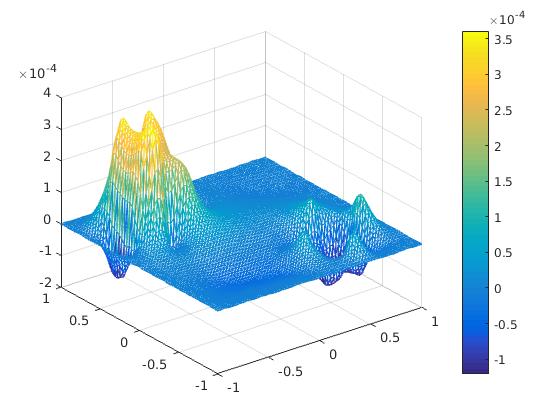}
\includegraphics[scale=0.2]{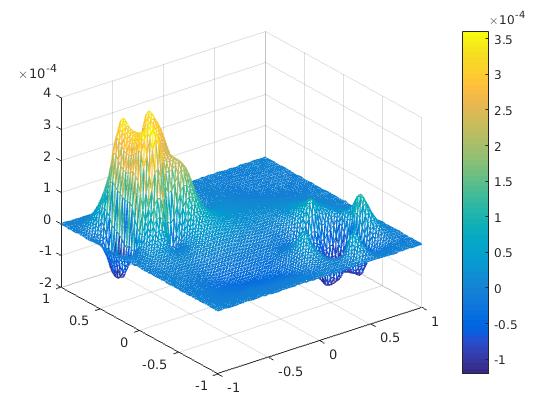}
\includegraphics[scale=0.2]{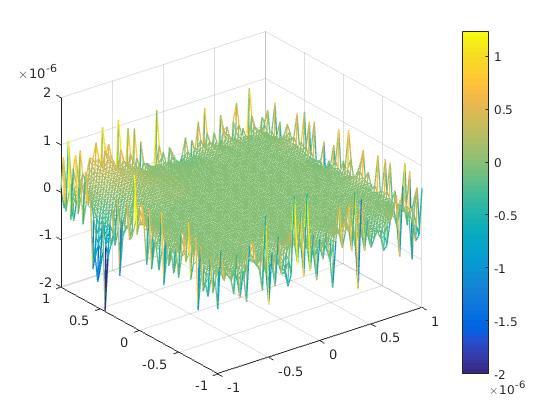}
\end{center}
\caption{Computed numerical solution $f_\ell$ of the algorithm at the final iteration, and the differences $\mathcal{N}^{h_\ell}_{ f^\dag}  j^\dag - \mathcal{N}^{h_\ell}_{f_\ell} j_{\delta_\ell}$, ~$ \mathcal{D}^{h_\ell}_{f^\dag} g^\dag-\mathcal{D}^{h_\ell}_{f_\ell} g_{\delta_\ell}$ and $\mathcal{D}^{h_\ell}_{f_\ell} g_{\delta_\ell} - \mathcal{N}^{h_\ell}_{f_\ell} j_{\delta_\ell}$.}
\label{h1}
\end{figure}
\end{example}

\begin{example}\label{ex_mult}
In present example we assume that multiple measurements are available, say $\left(j_\delta^i,g_\delta^i \right)_{i=1,\ldots,I}$. Then, problem $\left(\mathcal{P}^h_{\rho,\delta}\right)$ in Section \ref{Stability} is given by
$$\min_{f\in L^2(\Omega)} \bar{\Upsilon}^h_{\rho,\delta} (f)
:= \min_{f\in L^2(\Omega)} \left( \underbrace{\frac{1}{I}\sum_{i=1}^I\int_\Omega Q\nabla \left(\mathcal{N}^h_fj_\delta^i-
\mathcal{D}^h_fg_\delta^i\right) \cdot \nabla \left(\mathcal{N}^h_fj_\delta^i-
\mathcal{D}^h_fg_\delta^i\right)}_{:= \bar{\mathcal{J}}_\delta^h (q)} + \rho  \left \| f-f^*
\right \|^2_{L^2(\Omega)}\right)  \eqno\left(\bar{\mathcal{P}}^h_{\rho,\delta}\right),$$
which also attains a solution $\bar{f}^h_{\rho,\delta}$. The Neumann boundary condition in the equation \eqref{17-5-16ct3*} is chosen in the same form as \eqref{20-6-17ct5}, i.e.
\begin{equation}\label{20-6-17ct1}
\begin{aligned}
j^\dag_{(A,B,C,D)}
&= A\cdot \chi_{(0,1]\times\{-1\}} - A\cdot\chi_{[-1,0]\times\{1\}} + B\cdot\chi_{(0,1]\times\{1\}} -B\cdot\chi_{[-1,0]\times\{-1\}} \\
&~\quad +C\cdot\chi_{\{-1\}\times(-1,0]} - C\cdot\chi_{\{1\}\times(0,1)}  + D\cdot\chi_{\{1\}\times(-1,0]} - D\cdot\chi_{\{-1\}\times(0,1)},
\end{aligned}
\end{equation}
and depends on the constants $A,B,C$ and $D$. Let $g^\dag_{(A,B,C,D)} := \gamma\mathcal{N}_{f^\dag}j^\dag_{(A,B,C,D)}$ and assume that noisy observations are given by
\begin{align}\label{20-6-17ct2}
\left( j^{(A,B,C,D)}_{\delta_{\ell}}, g^{(A,B,C,D)}_{\delta_{\ell}} \right) = \left( j^\dag_{(A,B,C,D)}+ \theta\cdot R_{j^\dag_{(A,B,C,D)}}, g^\dag_{(A,B,C,D)}+ \theta\cdot R_{g^\dag_{(A,B,C,D)}}\right),
\end{align}
where $R_{j^\dag_{(A,B,C,D)}}$ and $R_{g^\dag_{(A,B,C,D)}}$ denote $\partial M^{h_\ell}\times 1$-matrices of random numbers on the interval $(-1,1)$. Different from \eqref{3-7-17ct1}, the constant $\theta$ appeared in the equation \eqref{20-6-17ct2} is now independent of the grid level $\ell$.

In the case $(A,B,C,D)=(1,2,3,4)$ we have a single noisy measurement couple, i.e. $I=1$. We now fix $D=4$, and let $(A,B,C)$ take all permutations $\mathcal{S}_3$ of the set $\{1,2,3\}$. Then, the equations \eqref{20-6-17ct1}--\eqref{20-6-17ct2} generate $I=6$ measurements. Similarly, if $(A,B,C,D)$ takes all permutations $\mathcal{S}_4$ of $\{1,2,3,4\}$ we get $I=16$ measurements. With $\theta=0.1$ and $\ell=64$ we compute the noise level
$$
\bar{\delta}_\ell= \begin{cases}
 \big\|j^{(1,2,3,4)}_{\delta_\ell} -j^\dag_{(1,2,3,4)}\big\|_{L^2(\partial\Omega)} + \big\|g^{(1,2,3,4)}_{\delta_\ell} -g^\dag_{(1,2,3,4)}\big\|_{L^2(\partial\Omega)} \quad\mbox{if}\quad (A,B,C,D)=(1,2,3,4),\\
\dfrac{1}{6}\sum_{(A,B,C)\in\mathcal{S}_3} \big\|j^{(A,B,C,4)}_{\delta_\ell} -j^\dag_{(A,B,C,4)}\big\|_{L^2(\partial\Omega)} + \big\|g^{(A,B,C,4)}_{\delta_\ell} -g^\dag_{(A,B,C,4)}\big\|_{L^2(\partial\Omega)} \quad\mbox{if}\quad D=4,\\
\dfrac{1}{16}\sum_{(A,B,C,D)\in\mathcal{S}_4} \big\|j^{(A,B,C,D)}_{\delta_\ell} -j^\dag_{(A,B,C,D)}\big\|_{L^2(\partial\Omega)} + \big\|g^{(A,B,C,D)}_{\delta_\ell} -g^\dag_{(A,B,C,D)}\big\|_{L^2(\partial\Omega)}.
\end{cases}$$
The corresponding numerical results for the multiple measurement case are presented in the Table \ref{b4}.
\begin{table}[H]
\begin{center}
\begin{tabular}{|c|l|l|l|l|l|l|l|l|}
\hline \multicolumn{9}{|c|}{ {\bf Numerical results for $\ell=64$, $\theta=0.1$ with multiple observations} }\\
\hline $I$  &\scriptsize {\bf Iterate} &\scriptsize {\bf Tolerance} &\scriptsize $\bar{\delta}_\ell$ &\scriptsize $L^2_f$ &\scriptsize $L^2_\mathcal{N}$ &\scriptsize $L^2_\mathcal{D}$ &\scriptsize $H^1_{\mathcal{N}}$ &\scriptsize $H^1_{\mathcal{D}}$\\
\hline
1  & 531& -3.2313e-8& 0.3292& 0.3280&  5.9096e-3 &5.9090e-3 &0.1225&0.1221\\
\hline
6 & 517& -7.1620e-9 & 0.3331&  0.2583& 4.3125e-3& 4.3122e-3 &7.9322e-2&7.9320e-2\\
\hline
16 & 536 & -6.4706e-8& 0.3289&  0.1747& 2.8465e-3& 2.8461e-3 &5.2318e-2&5.2314e-2\\
\hline
\end{tabular}
\caption{Numerical results for $\ell=64$, $\theta=0.1$, and with multiple measurements $I=1,6,16$.}
\label{b4}
\end{center}
\end{table}
Finally, in Figure \ref{h2} from left to right we show the interpolation $I_1^{h_{\ell}} f^\dag$ of the exact source and the computed numerical solution $q_\ell$ of the algorithm at the final iteration for $\ell=64$, $\theta=0.1$, and $I=16,6,1$, respectively.
\begin{figure}[H]
\begin{center}
\includegraphics[scale=0.2]{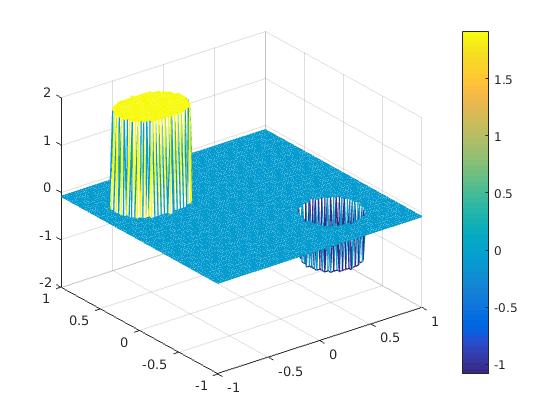}
\includegraphics[scale=0.2]{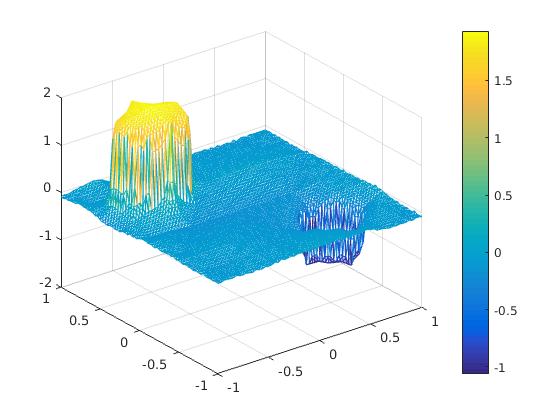}
\includegraphics[scale=0.2]{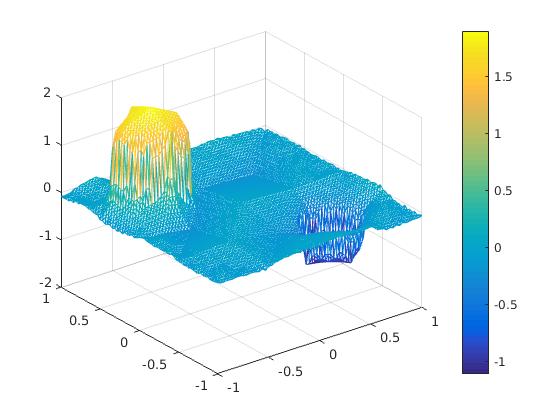}
\includegraphics[scale=0.2]{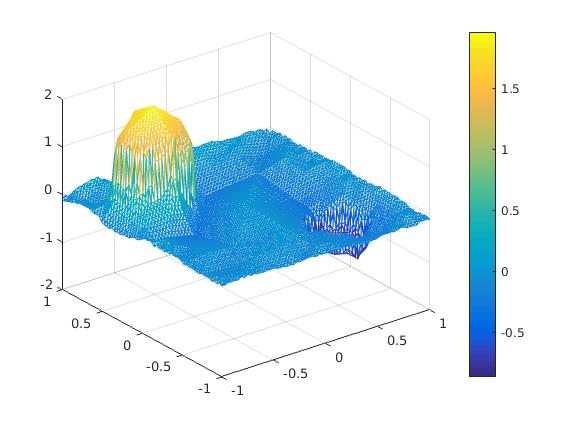}
\end{center}
\caption{Interpolation $I^{h_\ell}_1 f^\dag$, computed numerical solution $f_\ell$ of the algorithm at the final iteration for $\ell = 64$, $\theta = 0.1$, and with multiple measurements $I = 16, 6, 1$, respectively.}
\label{h2}
\end{figure}
\end{example}

\section*{Acknowledgments}

%The authors M. Hinze, B. Hofmann and T.N.T. Quyen would like to thank two anonymous referees and the editor for their constructive comments and suggestions which helped significantly to improve the present paper.

The authors thank the Referee and Editor for their valuable comments and suggestions.

M. Hinze gratefully acknowledges support of the Lothar Collatz Center for Computing in Science at the
University of Hamburg. \\
B. Hofmann gratefully acknowledges support by the German Research Foundation (DFG) under grant HO~1454/12-1. \\
T. N. T. Quyen gratefully acknowledges support of the Alexander von Humboldt Foundation.

\end{document}